\documentclass[12pt]{amsart}
\usepackage{caption,amsmath,amssymb,amsbsy,amsfonts,amsthm,latexsym,graphicx,multirow,color,hyperref,enumerate,srcltx}
\newtheorem{lem}{Lemma}[section]%
\newtheorem{theorem}[lem]{Theorem}%
\newtheorem{defi}[lem]{Definition}%
\newtheorem{exam}[lem]{Example}%
\newtheorem{prop}[lem]{Proposition}%

\newtheorem{lemma}{Lemma}[section]%
\newtheorem{definition}[lem]{Definition}%

\def\a{\alpha} \def\b{\beta}  \def\d{\delta} 
\def\r{\rho} \def\s{\sigma}

\def\G{{\it\Gamma}}
\def\D{{\it \Delta}}

\def\Om{{\it \Omega}}

  \def\ov{\overline V}

 \def\lg{\langle} \def\rg{\rangle}
\def\nd{\mathrel{\bigm|\kern-.7em/}}

\def\f{\noindent}

\def\PSL{\hbox{\rm PSL}}
\def\PSp{\hbox{\rm PSp}}
\def\PSU{\hbox{\rm PSU}}
\def\GO{\hbox{\rm GO}}

\def\SL{\hbox{\rm SL}}
\def\GL{\hbox{\rm GL}}
\def\PGL{\hbox{\rm PGL}}

\def\Aut{\hbox{\rm Aut\,}}
\def\Inn{\hbox{\rm Inn}}

\def\soc{\hbox{\rm soc}}

\def\Cay{\hbox{\rm Cay}}

\def\Out{\hbox{\rm Out}}
\def\D{\Delta}

\def\Ga{{\it\Gamma}}
\def\K{{\bf K}}

\def\calB{{\mathcal B}}

\def\Z{{\bf Z}}

\def\PL{{\rm P\Gamma L} }

\def\ov{\overline}
\def\l{\langle}
\def\r{\rangle}

\def\demo{{\bf Proof}\hskip10pt}

\def\mz{{\mathbb Z}}
\def\qed{\hskip10pt $\Box$\vspace{3mm}}

\begin{document}
\title[Homogeneous graphs]{Finite $3$-connected homogeneous graphs}
\author[Li]{Cai Heng Li}
\address{Cai Heng Li\\
Department of Mathematics\\
Southern University of Science and Technology\\
Shenzhen 518055, Guangdong, P.R. China}
\email{lich@sustc.edu.cn}

\author[Zhou]{Jin-Xin Zhou}
\address{Jin-Xin Zhou\\
School of Mathematics and Statistics\\
Beijing Jiaotong University\\
Beijing 100044, P.R. China}
\email{jxzhou@bjtu.edu.cn}
\date{}

\begin{abstract}
A finite graph $\G$ is said to be {\em $(G,3)$-$($connected$)$ homogeneous} if every isomorphism between any two isomorphic (connected) subgraphs of order at most $3$ extends to an automorphism $g\in G$ of the graph, where $G$ is a group of automorphisms of the graph.
In 1985, Cameron and Macpherson determined all finite $(G, 3)$-homogeneous graphs. In this paper, we develop a method for characterising  $(G,3)$-connected homogeneous graphs. It is shown that for a finite $(G,3)$-connected homogeneous graph $\G=(V, E)$, either $G_v^{\G(v)}$ is $2$-transitive or $G_v^{\G(v)}$ is of rank $3$ and $\G$ has girth $3$, and that the class of finite, non-multipartite, $(G,3)$-connected homogeneous graphs is closed under taking normal quotients. This leads us to study graphs where $G$ is quasiprimitive on $V$. We determine the possible quasiprimitive types for $G$ in this case and give new constructions of examples for some possible types.

\bigskip

\end{abstract}
\maketitle

\section{Introduction}


A graph is called {\em homogeneous} if any isomorphism between finite induced subgraphs extends to an automorphism of the graph. 
In 1976, Gardiner~\cite{Gardiner1976} gave an explicit classification of the finite homogeneous graphs, and later, Lachlan and Woodrow~\cite{Lachlan-Woodrow} extended this to the infinite countable homogeneous graphs. If we only consider the connected subgraphs, the connected homogeneity, a natural weakening of homogeneity, arises. A graph $\G$ is {\em connected homogeneous} if every isomorphism between connected induced subgraphs extends to an automorphism of $\G$. The connected homogeneous graphs lie between homogeneous graphs and distance-transitive graphs. The finite and infinite countable connected homogeneous graphs have also been classified by Gardiner et al.~\cite{Enomoto,Gardiner1978,Gray-Macpherson}. For more results regarding the (connected) homogeneous graphs, we refer the reader to \cite{HH,Gray,GMPR}.

Very few families of graphs are (connected) homogeneous. For example, a finite graph is homogeneous if and only if it is a union of isomorphic complete graphs, a regular complete multipartite graph, a grid of order 9, or a pentagon, see \cite{Gardiner1976}. A finite connected graph is connected homogeneous if and only if it is a homogeneous graph, a cycle of length greater than $5$, the Cartesian product $\K_n\square \K_n$ where $n\geq 3$, the strong product $\K_n\times\K_2$ where $n\geq 4$, the Petersen graph or the folded $5$-cube, see \cite[Theorem 3(a)]{Gardiner1978}. In view of these results, it is natural to relax the (connected) homogeneity. The following two symmetries have been considered in the literature.

\begin{definition}\label{Def-homog}
{\rm
Let $k$ be a positive integer.
Let $\Ga=(V,E)$ be a connected graph, and let $G\leq\Aut\Ga$.
\begin{itemize}
\item[(1)] If each isomorphism between any two isomorphic induced subgraphs of $\Ga$ of order at most $k$
extends to an automorphism $g$ of $\Ga$ such that $g\in G$, then $\Ga$ is called a {\it $(G,k)$-homogeneous} graph.


\item[(2)] If each isomorphism between any two isomorphic connected induced subgraphs of $\Ga$ of order at most $k$
extends to an automorphism $g$ of $\Ga$ such that $g\in G$, then $\Ga$ is called a {\it $(G,k)$-connected-homogeneous} graph, or simply called a {\it $(G,k)$-CH graph} for short.


\end{itemize}

}
\end{definition}

For $k=1$, each of these two types of symmetries is equivalent to the vertex-transitivity. When $k\geq 2$, finite $k$-homogeneous graphs have been classified in a collection of papers. For the case where $k=2$, a graph $\Ga$ is $(G, 2)$-homogeneous if and only if $\Ga$ and its complement $\ov\Ga$ are $G$-arc-transitive (namely, $G$ is transitive on the arcs of $\Ga$ and $\ov\Ga$). Thus either it is a complete multipartite graph $\K_{m[b]}$, or it is an orbital graph of a primitive permutation group
of rank 3. Due to the classification of primitive permutation groups of rank 3 given in \cite{Bannai,Kantor,Liebeck-affine,Liebeck-Saxl},
$(G, 2)$-homogeneous graphs are in some sense known. For the case where $k>2$, all finite $k$-homogeneous graphs are explicitly known. Actually, every finite $5$-homogeneous graph is homogeneous (see \cite{Cameron-6-transitive}), and the only finite $4$-homogeneous but not $5$-homogeneous graphs are the Schl$\ddot{{\rm a}}$fly graph and its complement (see \cite{Buczak} or Note added in the proof of \cite{Cameron-6-transitive}). A list of finite $3$-homogeneous graphs was given in \cite{Cameron-3-hom}.




Now let us consider $k$-CH graphs with $k\geq 2$. Recall that an {\em arc\/} in a graph is an ordered pair of adjacent vertices.
Similarly, a {\em $2$-arc\/} in a graph $\G$ is an ordered triple $(u, v, w)$ of three distinct vertices of $\G$ such that
$v$ is adjacent to both $u$ and $w$. By definition, a graph $\Ga$ is $2$-CH if and only if $\Ga$ is arc-transitive. There is no hope to obtain a classification of arc-transitive graphs because this class is too rich. It is natural to study $k$-CH graphs for $k\geq3$. Moreover, a graph of girth at least 4 is $3$-CH if and only if it is $2$-arc transitive. A fair amount of work have been done on the $2$-arc-transitive graphs in the literature, see, for example, \cite{Fang-Praeger-Suzuki,Fang-Praeger-Ree,GiudiciLCP,Ivanov-Praeger,Praeger1993-AJC,Praeger1993}. Motivated by this, we shall focus on 3-CH graphs of girth $3$. Notice that a {$2$-arc\/} $(u, v, w)$ in a graph $\G$ is called a $2$-{\em geodesic} if $u$ and $w$ are not adjacent. Clearly, every 3-CH graph of girth $3$ is $2$-geodesic transitive. The $2$-geodesic transitive graphs are also a class of interesting graphs and have been extensively studied, we refer the reader to \cite{DJLP-JCTA,DJLP-JACO}. In this paper, we aim to develop a method for characterising $(G,3)$-CH graphs.

It is well-known and easily shown that a vertex-transitive graph is $2$-arc-transitive if and only if it is locally $2$-transitive. The following proposition shows that the class of $3$-CH graphs, in some sense, corresponds to the locally rank $3$ action.

\begin{prop}\label{prop-local-action}
A connected graph $\Ga$ is $(G,3)$-CH for some $G\leq\Aut\G$ if and only if $G$ is vertex-transitive on $\Ga$,
and for a vertex $u$ of $\Ga$, either $G_u^{\Ga(u)}$ is $2$-transitive, or $G_u^{\Ga(u)}$ is of rank $3$ and the girth of $\Ga$ is $3$.
\end{prop}

Our next result is about the local structure of the $(G,k)$-CH graphs with $k\geq 3$.

\begin{theorem}\label{local-pty}
Let $\G=(V, E)$ be a connected $(G,k)$-CH graph for some $G\leq\Aut\G$ and $k\geq 3$. Then 
one of the following holds:
\begin{itemize}
\item[{\rm(i)}] $\Ga$ is isomorphic to the complete graph $\K_n$, and $G$ is $k$-transitive on $V$, where $n=|V|\geq k$.

\item[{\rm(ii)}] $\Ga$ is isomorphic to the complete multipartite graph $\K_{m[b]}$ with $m\geq 3$ parts of size $b\geq 2$, and $G\leq X\wr Y$, where $X$ is a $r_1$-transitive permutation group of degree $b$ with $r_1={\rm min}\{b, k-1\}$, $Y$ is a $r_2$-transitive permutation group of degree $m$ with $r_2={\rm min}\{m, k\}$, and $mb=|V|$.

\item[{\rm(iii)}] $\Ga$ is $(G,2)$-arc transitive of girth at least $4$.

\item[{\rm(iv)}] $[\Ga(u)]\cong r\K_b$ with $r,b>1$, and $G_u^{\Ga(u)}$ is imprimitive of rank $3$.

\item[{\rm(v)}] $G_u^{\Ga(u)}$ is primitive of rank $3$, and $[\Ga(u)]$ is connected. Furthermore, if $k\geq 4$ then the graphs are determined in \cite[Theorem~1.2]{DFLPZ}.
\end{itemize}
\end{theorem}

\f{\bf Remark~1}\ The graphs in (v) for which $G_u^{\Ga(u)}$ is transitive on the $2$-subsets of $\Ga(u)$  have been determined in \cite{J-locally-tria}. \medskip

We would like to propose the following problem.\medskip

\f{\bf Problem}\ Classify the graphs in (v) of Theorem~\ref{local-pty} for $k=3$.\medskip

One of the basic strategies in the study of symmetry in graphs is to study the normal quotients of graphs.
Let $\G$ be a vertex-transitive graph, and let $G\leq \Aut(\G)$ be vertex-transitive on $\G$.
Let $N$ be a non-trivial normal subgroup of $G$ which is intransitive on $V\G$.
The {\em quotient graph} $\G_N$ is defined as the graph with vertices being the $N$-orbits on $V\G$ such that
any two different vertices $B,C\in V\G_N$ are adjacent  if and only if there exist $u\in B$ and $v\in C$ which are adjacent in $\G$.
The original graph $\G$ is said to be a {\em cover} of $\G_N$ if $\G$ and $\G_N$ have the same valency.

The method of taking normal quotients has been very successful in investigating various families of graphs, for example,
$s$-arc transitive graphs~\cite{Praeger1993-AJC,Praeger1993} and locally $s$-arc-transitive graphs~\cite{GiudiciLCP},
where $s\ge2$. Finite 2-arc-transitive graphs form a subclass of 3-CH graphs, which is closed under taking normal quotients.
The following theorem shows that the class of finite $3$-CH graphs of girth $3$ is also closed under taking normal quotients.

\begin{theorem}\label{normal quotient}
Suppose that $\G$ is a connected $(G, 3)$-CH graph, which is not complete multipartite.
Let $N$ be a normal subgroup of $G$ which is intransitive on $V\G$.
Then either $\G_N=\K_2$ and $\G$ is $(G,2)$-arc-transitive, or $\G_N$ is $(G/N,3)$-CH and  $\G$ is a normal cover of $\G_N$.

\end{theorem}

\f{\bf Remark~2}\ Theorem~\ref{normal quotient} is not true for $(G, k)$-CH graph with $k\geq 4$. For example, it is easy to show that the $4$-cube $Q_4$ is $4$-CH, and the halved $4$-cube $\frac{1}{2}Q_4$ is not $4$-CH. However, $Q_4$ is a cover of $\frac{1}{2}Q_4$. \medskip

A transitive permutation group $G$ on a set $\Om$ is said to be {\em quasiprimitive}
if every nontrivial normal subgroup of $G$ acts transitively on $\Om$.  Praeger~\cite{Praeger1993} generalised the O'Nan-Scott theorem for primitive groups to quasiprimitive groups and showed that a finite quasiprimitive group is one of the following eight types. We shall describe these eight types along the lines of that in~\cite{Praeger}.

Let $G$ be a quasiprimitive permutation group on a finite set $\Om$ and take $\a\in\Om$. The {\em socle} of $G$, denoted by $\soc(G)$, is
the product of all minimal normal subgroups of $G$. Then $G$ has at most two minimal normal subgroups, and $\soc(G)\cong T_1\times\cdots\times T_d =T^d$, where $T_1\cong\cdots\cong T_d\cong T$ with $T$ a simple group. \medskip

\f{HA ({\em holomorph affine})}:\ {$\soc(G)$ is an abelian minimal normal subgroup, and $T=\mz_p$ for some prime $p$, and
$$N\lhd G\leq N: \Aut(N)\cong \mz_p^d: \GL(d,p)={\rm AGL}(d,p).$$
In this case, $G_\a$ is an irreducible subgroup of $\GL(d, p)$. }\medskip

\f{HS ({\em holomorph simple})}:\ $\soc(G)=M\times N$, where $M\cong N\cong T$ and $M,N$ are minimal normal subgroups. In this case,
we have $T\times T\leq G\leq T: \Aut(T)$. \medskip

\f{HC ({\em holomorph compound})}:\ $\soc(G)=M\times N$, where $M\cong N\cong T^\ell (\ell>1)$ and $M,N$ are minimal normal subgroups. In this case, we have $N\times N\leq G\leq N: \Aut(N)$.

\f{AS ({\em almost simple})}:\ $\soc(G)=T$ is nonabelian simple and $T\leq G\leq\Aut(T)$.\medskip

The remaining case is where $\soc(G)=T_1\times\cdots\times T_d=T^d (d>1)$ is a nonabelian minimal normal subgroup of $G$. There are four different types of primitive permutation groups. \medskip

\f{PA ({\em product action})}:\  $N=\soc(G)$ has no normal subgroup which is regular on $\Om$.\medskip

\f{TW ({\em twisted wreath product})}:\ $N=\soc(G)$ is minimal normal in $G$ and regular on $\Om$.  \medskip

\f{SD ({\em simple diagonal})}:\ $N=\soc(G)$ is not regular on $\Om$ and has a normal subgroup which is regular on $\Om$, and $N_\a\cong T$.

\f{CD\ ({\em compound diagonal})}:\ $N=\soc(G)$ is not regular on $\Om$ and has a normal subgroup $L$ which is regular on $\Om$. If $L\cong T^m (m\leq d-2)$, then $N_\a\cong T^{d-m}$.
\medskip

In Theorem~\ref{normal quotient}, if we choose $N$ to be maximal by inclusion subject to being intransitive on $V\G$,
then $G/N$ is quasiprimitive on $V\G_N$ in addition to $\G_N$ being $(G/N,3)$-CH.
Therefore, it is natural to consider the $(G,3)$-CH graphs such that $G$ is quasiprimitive on $V\G$.
The following theorem determines the possibilities for the type of $G$.

\begin{theorem}\label{th-quasiprimitive}
Suppose that $\G=(V,E)$ is a $(G,3)$-CH graph which is neither complete nor complete multipartite, and $G$ is quasiprimitive on $V$.
Then the following hold.
\begin{enumerate}[{\rm (1)}]
  \item $G$ is not of type {\rm HS, HC} or {\rm CD}.
  \item If $G$ is of type {\rm SD}, then $\G\cong\G_{A_5^3}, \G_{A_5^4}, \G_{U(3,4)^3}$ or $\G_{U(3,4)^4}$ $($see {\rm Examples~\ref{exam:sd-a5-1}--\ref{exam:sd-u34-2})}.
  \item There do exist Examples of $(G,3)$-CH graphs in case $G$ is of type {\rm HA, AS, PA} or {\rm TW}. $($see {\rm Section~\ref{sec:examples})}.
\end{enumerate}
\end{theorem}

The rest of the paper is organised as follows. In Section~\ref{sec:prelimiaries}, we give some notation and definitions. In Section~\ref{sec:locall-structure}, we prove Proposition~\ref{prop-local-action} and Theorem~\ref{local-pty}; and in Section~\ref{sec:quotient}, we prove Theorem~\ref{normal quotient}. In Section~\ref{sec:normalCayley}, we give a characterisation of $3$-CH normal Cayley graphs which will be used in the proof of Theorem~\ref{th-quasiprimitive}~(1) and (3), and in Section~\ref{sec:2-transitive}, we provide some properties of $2$-transitive groups which will be used in the proof Theorem~\ref{th-quasiprimitive}~(1) and (2). In Section~\ref{sec:part1}, we prove Theorem~\ref{th-quasiprimitive}~(1), in Section~\ref{sec:part2}, we prove Theorem~\ref{th-quasiprimitive}~(2), and in Section~\ref{sec:examples}, we prove Theorem~\ref{th-quasiprimitive} (3) by giving some examples of finite $(G,3)$-CH graphs such that $G$ is quasiprimitive on the vertices of type HA, AS, PA or TW. Finally, in Section~\ref{sec:proof}, we give the proof of Theorem~\ref{th-quasiprimitive}.

\section{Preliminaries}\label{sec:prelimiaries}

All groups considered in this paper are finite, and all graphs are finite, connected, simple and
undirected, unless explicitly stated. For the group-theoretic and graph-theoretic terminology not defined here we refer the reader to \cite{BMBook,WI}.

For a positive integer $n$, the expression $\mz_n$ denotes the cyclic group of order $n$, $D_{2n}$ denotes the dihedral group of order $2n$,
$A_n, S_n$ denote the alternating group and symmetric group of degree $n$, respectively. Before proceeding, we introduce some notation. For a positive integer $n$, $\K_n$ is the complete graph on $n$ vertices. The complete multipartite graph $\K_{m[b]}$ with $m\geq3,b\geq2$, has vertex set consisting of $m$ parts of size $b$, with edges between all pairs of vertices from distinct parts.

Let $\G=(V, E)$ be a graph with vertex set $V$ and edge set $E$.
An {\em arc} is an ordered pair of adjacent vertices, and an edge $\{u,v\}$ corresponds to arcs $(u,v)$ and $(v,u)$.
For any subset $B$ of $V$, the subgraph of $\G$ induced on $B$ will be denoted by $\G[B]$, and
when no confusion arises, it is simply written as $[B]$.
If $u,v\in V$, then $d(u,v)$ denotes the distance between $u$ and $v$ in $\G$.
The {\em diameter} $d$ of $\G$ is the maximal distance between two vertices in $\G$.
We shall assume that $d\geq1$.
For a vertex $v$, we write $\G_i(v)=\{u\ |\ d(u,v)=i\}$ for $1\le i\le d$;
and $\Ga_1(v)$ is simply denoted by $\Ga(v)$, which is the neighbourhood of $v$.

Let $G$ be a permutation group on a finite set $\Om$.
Let $\D$ be a subset of $\Om$.
Denote by $G_{\D}$ and $G_{(\D)}$ the subgroups of $G$ fixing $\D$ setwise and pointwise, respectively.
Let $G^{\D}$ represent the permutation group on $\D$ induced by $G$.
Then $G^{\D}\cong G_{\D}/G_{(\D)}$.
The group $G$ is {\em semiregular} on $\Om$  if the only element fixing a point in $\Om$ is the identity element of $G$.
We say that $G$ is {\em regular} on $\Om$ if it is both transitive and semiregular on $\Om$.

Let $G$ be a finite group and let $g\in G$. Let $\rho_g$ and $\s_g$ be the permutations of $G$ defined by $\rho_g: x\mapsto xg$ and $\s_g: x\mapsto g^{-1}xg$ for $x\in G$. The {\em right regular representation} of $G$ is the subgroup of Sym$(G)$ defined by $G_R:=\{\rho_g\ |\ g\in G\}$.
The map $\s_g$ is called an {\em inner automorphism} of $G$ induced by $g$, and the set of all $\s_g$, denoted by $\Inn(G)$, is called the {\em inner automorphism group} of $G$. Denote by $\Aut(G)$ the automorphism group of $G$. 

Given a finite group $H$ and a  subset $S$ of $H$ which does not contain the identity 1 of $H$,
the {\em Cayley graph} $\Cay(H,S)$ on $H$ with respect to $S$ is defined to be the graph with vertex set $H$
and edge set $\{\{g,sg\}\mid g\in H,s\in S\}$.
By the definition, we have the following basic properties, where $\Ga=\Cay(H,S)$:

\begin{itemize}
\item[(i)] $\Ga$ is undirected if and only if $S$ is self-inverse, namely, $g\in S$ if and only if $g^{-1}\in S$;

\item[(ii)] a Cayley graph $\Ga$ is connected if and only if $S$ generates $H$;

\item[(iii)] the right multiplications of elements of $H$ form a subgroup of  $\Aut\Ga$
which is vertex-transitive on $\Ga$; in particular, Cayley graphs are vertex-transitive.
\end{itemize}

There is a criterion to determine a graph to be a Cayley.

\begin{lemma}\label{Cay-criterion} {\rm(\cite[Proposition 16.3]{B})}
A graph $\Ga=(V, E)$ is a Cayley graph of a group $H$ if and only if $\Aut\Ga$ contains a subgroup which is isomorphic to $H$
and regular on $V$.
\end{lemma}

For a Cayley graph $\Ga=\Cay(H,S)$, let
\[\Aut(H,S)=\{\a\in\Aut(H)\ |\ S^\a=S\}.\]
It is easily shown that each element of $\Aut(H,S)$ induces an automorphism of $\Ga$, and  normalises $H_R$.
Moreover, we have the following lemma due to Godsil.

\begin{lemma}\label{Aut(H,S)}{\rm(\cite{Godsil,X1})}
The normaliser $N_A(H_R)=H_R: \Aut(H,S)$.
\end{lemma}

Let $G$ be a group acting transitively on a finite set $\Om$. A nonempty subset $\D$ of $\Om$ is called a {\em block} for $G$ if for each $g\in G$, either $\D^g=\D$ or $\D^g\cap \D=\emptyset$. We call $\{\D^g\ |\ g\in G\}$ a {\em block system} for $G$. Then $G$ is said to be {\em primitive} if the only blocks for $G$ are the singleton subsets or the whole of $\Om$. It is well-known that the orbits of a normal subgroup of $G$ form block system for $G$. If each nontrivial normal subgroup of $G$ is transitive on $\Om$, then $G$ is said to be {\em quasiprimitive}.

We now introduce a characterisation of primitivity. Consider the natural action of $G$ on the cartesian product $\Om\times\Om$. The orbits of $G$ on this set are called the {\em orbitals} of $G$ on $\Om$. The orbital $\D_1=\{(u,u)\ |\ u\in\Om\}$ is called the {\em diagonal orbital} and all other orbitals are said to be {\em non-trivial}. For each orbital $\D$, there is an orbital, denoted by $\D^*$, so that $(u,v)\in\D^*$ if and only if $(v,u)\in\D$. An orbital is {\em self-paired} if $\D^*=\D$. For each orbital $\D$ of $G$, the digraph Graph$(\D)$ is a digraph with vertex set $\Om$ and edge set $\D$. It is easy to see that $\D$ is self-paired if and only if Graph$(\D)$ is a graph. By \cite[Theorem~3.2A]{Dixon}, $G$ is primitive on $\Om$ if and only if Graph$(\Delta)$ is connected for each non-trivial orbital $\D$. For each orbital $\D$ of $G$ and each $u\in\Om$, define $\D(u)=\{v\in\Om\ |\ (u,v)\in\D\}$. Then the mapping $\D\mapsto\D(u)$ is a bijection from the set of orbitals of $G$ onto the set of orbits of $G_u$. In particular, the number of orbitals of $G$ is equal to the number of orbits of $G_u$; this number is called the {\em rank} of $G$. An orbit of $G_u$ for any $u\in\Om$ is called a {\em suborbit} of $G$, and if $\D$ and $\D^*$ are { paired } orbitals, then $\D(u)$ and $\D^*(u)$ are called {\em paired suborbits.}

\section{Proof of Proposition~\ref{prop-local-action} and Theorem~\ref{local-pty}}\label{sec:locall-structure}

In this section, we study the local structure of $k$-CH graphs with $k\geq 3$.
In particular, we shall prove Proposition~\ref{prop-local-action} and Theorem~\ref{local-pty}.

\vskip0.1in

\f{\bf Proof of Proposition~\ref{prop-local-action}:}\ Let $\G=(V,E)$ be a connected graph, let $G\leq\Aut\G$ and let $u$ be a vertex of $\G$.

Suppose that $\G$ is $(G,3)$-CH. Then $G$ is transitive on the vertices and the arcs of $\G$. First assume that $[\G(u)]$ is either a complete graph or an empty graph. For any two pairs $(v, w)$ and $(v', w')$ of neighbours of $u$, either both $\{v, w\}$ and $\{v', w'\}$ are edges of $\G$ or both $\{v, w\}$ and $\{v', w'\}$ are non-edges of $\G$. Let $\Sigma$ and $\Sigma'$ be the subgraphs of $\G$ induced by $\{v, u, w\}$ and $\{v', u, w'\}$, respectively. Then $\Sigma$ and $\Sigma'$ are both isomorphic to a triangle or a 2-path. In particular, there exists an isomorphism, say $\s$, between $\Sigma$ and $\Sigma'$ such that $(v, u, w)^\s=(v', u, w')$. Since $\G$ is $(G,3)$-CH, there exists $g\in G$ sending $(v, u, w)$ to $(v', u, w')$, and then $g\in G_u$ and $g$ sends $(v, w)$ to $(v', w')$. It follows that $G_u^{\G(u)}$ is $2$-transitive. Next assume that $[\G(u)]$ is neither a complete graph nor an empty graph. Then $\G$ is a non-complete graph of girth $3$. Since $G$ is transitive on the arcs of $\G$, $G_u^{\G(u)}$ is transitive. Let $(v_1, w_1)$ and $(v_2, w_2)$ be two arcs of the subgraph $[\Ga(u)]$. Then $\{v_1, u, w_1\}$ and $\{v_2, u, w_2\}$ induces two triangles, and hence there exists an isomorphism, say $\s$, between these two triangles such that $(v, u, w)^\s=(v', u, w')$. Since $\G$ is $(G,3)$-CH, there exists an automorphism $g\in G$ such that $(v_1, u, w_1)^g=(v_2, u, w_2).$
Then $g\in G_u$ and maps $(v_1, w_1)$ to $(v_2, w_2)$, and so $G_u$ is arc-transitive on $[\Ga(u)]$.
Assume now that $(v_1, w_1)$ and $(v_2, w_2)$ are two non-adjacent pairs of vertices of the subgraph $[\Ga(u)]$.
Then $(v_1, u, w_1)$ and $(v_2, u, w_2)$ form two 2-geodesics, and similarly, $G_u$ is transitive on the set of non-adjacent pairs of vertices of $[\Ga(u)]$. Thus $G_u^{\G(u)}$ has exactly two non-trivial orbitals on $\G(u)$, and
hence $G_u^{\G(u)}$ is of rank $3$.

Conversely, suppose that $G$ is vertex-transitive on $\G$, and that either $G_u^{\Ga(u)}$ is $2$-transitive, or $G_u^{\Ga(u)}$ is of rank $3$ and the girth of $\Ga$ is $3$. Then $\G$ is also $G$-arc-transitive.

Let $\Sigma_1$ and $\Sigma_2$ be two induced subgraphs of $\G$ of order 3 which are connected and isomorphic,
and let $\s$ be an isomorphism between $\Sigma_1$ and $\Sigma_2$.
Then $\Sigma_i$ is a 2-geodesic or a triangle, where $i=1$ or $2$, and hence $\Sigma_i$ contains a 2-arc $(v_i, u_i, w_i)$
such that $(v_1,u_1,w_1)^\s=(v_2,u_2,w_2)$. Since $G$ is transitive on the vertices, there exists $g_i\in G$ such that $u_i^{g_i}=u$ with $i=1,2$. So $(v_1, u_1, w_1)^{g_1}=(v_1^{g_1}, u, w_1^{g_1})$ and $(v_2, u_2, w_2)^{g_2}=(v_2^{g_2}, u, w_2^{g_2})$.

If $G_u^{\Ga(u)}$ is $2$-transitive, then there exists $g\in G_u$ mapping $(v_1^{g_1}, w_1^{g_1})$ to $(v_2^{g_2}, w_2^{g_2})$. So $(v_1,u_1,w_1)^{g_1gg_2^{-1}}=(v_2,u_2,w_2)=(v_1,u_1,w_1)^\s.$ So $\G$ is $(G,3)$-CH.

If $G_u^{\Ga(u)}$ is of rank $3$ and the girth of $\Ga$ is $3$, then $\G$ is $G$-arc-transitive, and
$G_u$ is transitive on the set of arcs, and also on the set of non-adjacent pairs of vertices in the subgraph $[\Ga(u)]$.
Since $\Sigma_1\cong\Sigma_2$, one has $\Sigma_1^{g_1}\cong\Sigma_2^{g_2}$, and hence $(v_1^{g_1}, w_1^{g_1})$ is an arc of $[\Ga(u)]$ if and only if $(v_2^{g_2}, w_2^{g_2})$ is an arc of $[\Ga(u)]$. It follows that there exists an element $g\in G_u$ such that $(v_1^{g_1}, w_1^{g_1})^g=(v_2^{g_2},w_2^{g_2})$. Again, we have $(v_1,u_1,w_1)^{g_1gg_2^{-1}}=(v_2,u_2,w_2)=(v_1,u_1,w_1)^\s,$ and hence $\G$ is $(G,3)$-CH.
\hfill\qed

We are ready to prove Theorem~\ref{local-pty}.
\vskip0.1in

\f{\bf Proof of Theorem~\ref{local-pty}:}\  Let $\Ga=(V,E)$ be a connected $(G, k)$-CH graph with $k\geq 3$.
If $\Ga\cong\K_n$ with $n=|V|$, then $n\geq k$ and $G$ is $k$-transitive on $V$, and so part (i) happens.

Assume that $\Ga\cong\K_{m[b]}$ with $m\geq 3$ and $b\geq 2$. Let $B_1, B_2, \ldots, B_m$ be the $m$ parts of $V(\Ga)$. Then $\mathcal{B}=\{B_1, B_2, \ldots, B_m\}$ is a block system of $G$ on $V(\G)$. Let $X=G_{B_1}^{B_1}$ and let $Y$ be the permutation group of $G$ induced on $\mathcal{B}$. Then $G\leq X\wr Y$.

Let $r_1={\rm min}\{b, k-1\}$. Take two arbitrary $r_1$-tuples of vertices in $B_1$, say $(v_1, v_2,\ldots, v_{r_1})$ and $(u_1, u_2,\ldots, u_{r_1})$. For any $w\in B_2$, let $\Sigma_1,\Sigma_2$ be the subgraphs induced by $\{w\}\cup\{v_1, v_2,\ldots, v_{r_1}\}$ and $\{w\}\cup\{u_1, u_2,\ldots, u_{r_1}\}$, respectively. Then $\Sigma_1\cong\Sigma_2\cong\K_{r_1,1}$, and there exists an isomorphism, say $\s$, from $\Sigma_1$ to $\Sigma_2$ such that $w^\s=w$ and $v_i^\s=u_i$ with $i=1,2,\ldots,r_1$. Clearly, $\Sigma_i$ has order $r_1+1\leq k$. Since $\Ga$ is $(G,k)$-CH, there exits $g\in G$ such that $w^g=w$ and $v_i^g=u_i$ with $i=1,2,\ldots,r_1$. Since $B_1$ is a block of imprimitivity of $G$ on $V(\Ga)$, one has $g\in G_{B_1}$ and hence $g^{B_1}\in X$ sends $(v_1, v_2,\ldots, v_{r_1})$ to $(u_1, u_2,\ldots, u_{r_1})$. So $X$ is $r_1$-transitive on $B_1$.

Now let $r_2={\rm min}\{m, k\}$. Take two arbitrary $r_2$-tuples of blocks in $\mathcal{B}$, say $(B_{i_1}, B_{i_2}, \cdots, B_{i_{r_2}})$ and $(B_{j_1}, B_{j_2}, \cdots, B_{j_{r_2}})$. Take $u_{i_t}\in B_{i_t}$ and $v_{j_t}\in B_{j_t}$ with $t=1,2,\ldots,r_2$. Let $\D_1,\D_2$ be the subgraphs induced by $\{u_{i_1}, u_{i_2}, \ldots, u_{i_{r_2}}\}$ and $\{v_{j_1}, v_{j_2}, \ldots, v_{j_{r_2}}\}$, respectively. Then $\D_1\cong\D_2\cong\K_{r_2}$, and there exists an isomorphism, say $\d$, from $\D_1$ to $\D_2$ such that $u_{i_t}^\s=v_{i_t}$ with $t=1,2,\ldots,r_2$. Clearly, $\D_i$ has order $r_2\leq k$. Since $\Ga$ is $(G,k)$-CH, there exists $g\in G$ such that $u_{i_t}^g=v_{i_t}$ with $t=1,2,\ldots,r_2$. Since $\mathcal{B}$ is a block system of $G$ on $V(\G)$, one has $(B_{i_1}, B_{i_2}, \cdots, B_{i_{r_2}})^g=(B_{j_1}, B_{j_2}, \cdots, B_{j_{r_2}})$. Thus, $Y$ is $r_2$-transitive on $\mathcal{B}$. This proves part (ii).

If $\Ga$ has girth at least $4$, then since $\G$ is $(G,k)$-CH with $k\geq 3$, by Proposition~\ref{prop-local-action}, $G$ is vertex-transitive on $\G$ and $G_u^{\G(u)}$ is $2$-transitive. It follows that $\G$ is $(G,2)$-arc-transitive. Thus, part (iii) happens.

In what follows, assume that $\G$ is neither a complete graph nor a complete multipartite graph, and assume further that $\Ga$ is of girth 3.
Let $u\in V$ be a vertex. Obviously, $[\Ga(u)]$ is $(G, k-1)$-homogeneous.
By Proposition~\ref{prop-local-action}, the permutation group $G_u^{\Ga(u)}$ is of rank 3.
Thus, all arcs of $[\Ga(u)]$ are equivalent under $G_u$, and
all pairs of non-adjacent vertices are equivalent under $G_u$.

First assume that the induced subgraph $[\Ga(u)]$ is disconnected.
Let $C_1,\dots,C_r$ be the connected components.
Suppose that $C_1$ is not a complete graph.
Then there exist two vertices $u,v$ of $C_1$ which are not adjacent.
Let $w$ be a vertex of $C_2$.
Since $G_u^{\Ga(u)}$ is of rank 3, $(u,v)$ and $(u,w)$ are equivalent under $G_u$, which is a contradiction.
Therefore, $C_1$ is a complete graph, and so is each $C_i$.
Hence $[\Ga(u)]=r\K_\ell$, as in part~(iv).

Next assume that $[\Ga(u)]$ is connected.
Suppose that $G_u^{\Ga(u)}$ is imprimitive.
Let $\calB=\{B_1,B_2,\dots,B_r\}$ be a non-trivial block system of $G_u^{\Ga(u)}$.
Then there is no edge lying inside a block $B_i$.
It implies that for any two distinct blocks $B_i$ and $B_j$, the induced subgraph $[B_i\cup B_j]$ is isomorphic to $\K_{\ell,\ell}$,
a complete bipartite graph, where $\ell=|B_i|$.
Thus the graph $[\Ga(u)]=\K_{r[\ell]}$, a complete multipartite graph, and the valency of $\Ga$ equals $r\ell$.
For any $v\in\G(u)$, considering the neighbourhood $\Ga(v)$ shows that the induced subgraph $[\Ga(v)]=\K_{r[\ell]}$.
It follows from \cite[Lemma~7]{Gardiner1978} that $\Ga=\K_{(r+1)[\ell]}$, which is a contradiction.
Therefore, $G_u^{\Ga(u)}$ is primitive. In this case, if $k\geq 4$, then all $(G,k)$-CH graphs have been determined in \cite[Theorem~1.2]{DFLPZ}. This proves part~(v) of Theorem~\ref{local-pty}.\hfill\qed

\section{Proof of Theorem~\ref{normal quotient}}\label{sec:quotient}

The class of 2-arc-transitive graphs is closed under taking normal quotient by \cite{Praeger1993}.
In this section we prove that the class of 3-CH graphs is closed under taking normal quotient,
which proves Theorem~\ref{normal quotient}. We first prove the following lemma.

\begin{lemma}\label{quotient}
Let $\Ga=(V,E)$ be a connected $(G,3)$-CH graph which is not a complete multipartite graph.
Let $N\lhd G$ have at least three orbits on $V$.
Then
\begin{itemize}
\item[{\rm (i)}] $\G$ is a normal cover of $\G_N$ and $N$ is semiregular on $V$, and

\item[{\rm (ii)}] $G_v$ acting on $\Ga(v)$ is permutationally isomorphic to
$\ov G_{\ov v}^{\Ga_N(\ov v)}$, where $\ov G=G/N$, and $\ov v$ is the vertex of $\Ga_N$ corresponding to $v$, and

\item[{\rm (iii)}] $\G_N$ is $(G/N, 3)$-CH.
\end{itemize}
\end{lemma}

\f\demo
Let $V_N$ be the set of $N$-orbits on $V$, and let $K$ be the kernel of $G$ acting on $V_N$.
Then $K\lhd G$,  $G/K$ is a subgroup of $\Aut\Ga_N$, and  $\Ga_N=\Ga_K$ is $G/K$-arc-transitive.
The valency of $\G_N$ is larger than 1. 
Since $\Ga$ is $G$-arc-transitive, it implies that there is no edge lying inside a block $v^N$, where $v\in V$.

Let $\{u,v\}$ be an edge of $\Ga$.
Since the orbit $v^N$ is a block of $G$, the intersection $v^N\cap\Ga(u)$ is a block of $G_u$ acting on $\Ga(u)$.
Suppose that $|v^N\cap\Ga(u)|\ge2$.
Let $w\in (v^N\cap\G(u))\setminus\{v\}$.
Since no edge lies in the same orbit of $N$, the vertices $v$ and $w$ are not adjacent, and $(v,u,w)$ is a 2-geodesic.
Since $\Ga_N$ has valency at least 2, there exist vertices which are in $\Ga(u)$ but not in $v^N$.
Let $w'\in\Ga(u)\setminus v^N$.
If $v,w'$ are not adjacent, then $(v,u,w')$ is a 2-geodesic, and so there is $g\in G$ such that
$(v,u,w)^g=(v,u,w')$ as $\Ga$ is $(G,3)$-CH, which contradicts that $v^N\cap \Ga(u)$ is a block of $G_u$.
Hence $v,w'$ are adjacent, and so $v$ is adjacent to all vertices of $\Ga(u)\setminus v^N$.
It implies that the induced subgraph $[\Ga(u)]$ is a complete multipartite graph. Since $G$ is vertex-transitive on $\Ga$, we have $[\G(u)]\cong[\G(v)]$ for all $v\in V$. It then follows from \cite[Lemma~7]{Gardiner1978} that $\Ga$ itself is a complete multipartite graph, which is a contradiction.
Therefore, $|v^N\cap\Ga(u)|=1$, namely, the orbit $v^N$ intersects the neighbourhood $\Ga(u)$
at the single vertex $v$.
So the quotient graph $\G_N$ and the original graph $\G$ have the same valency, and
hence $\G$ is a normal cover of $\G_N$.
In particular, $K=N$ is semiregular on $V$, and $N_v=1$. This proves part (i).

As mentioned above, the factor group $G/N$ is arc-transitive on the quotient graph $\Ga_N$.
For convenience, write $\ov G=G/N$, $\ov g=gN\in\ov G$ for any element $g\in G$, $\Sigma=\Ga_N$, and
for $v\in V$, the vertex $v^N$ of $\Sigma$ is denoted by $\ov v$.
Then the neighbourhood $\Sigma(\ov v)=\{\ov w \mid w\in\Ga(v)\}$, and the stabiliser $\ov G_{\ov v}=G_vN/N\cong G_v$.
Let $\phi$ be the isomorphism between $G_v$ and $\ov G_{\ov v}$ defined by
\[\phi:\ \ g\mapsto \ov g,\ \ \mbox{where $g\in G_v$}.\]
Label vertices of $\Ga(v)$ as $\{w_1,w_2,\dots,w_k\}$.
Then for any $w_i,w_j\in\Ga(v)$ and any $g\in G_v$, we have
\[w_i^g=w_j \Longleftrightarrow \ov w_i^{\ov g} =(w_i^N)^{gN}
  =w_i^{gg^{-1}NgN}=(w_i^g)^{g^{-1}NgN}=w_j^N=\ov w_j.\]
Hence $\ov G_{\ov v}$ acting on $\Sigma(\ov v)$ is permutationally isomorphic to $G_v$ acting on $\Ga(v)$. This proves part (ii).

Since $\G$ is $(G,3)$-CH, by Proposition~\ref{prop-local-action}, $G$ is vertex-transitive on $\G$ and either $G_u^{\Ga(u)}$ is $2$-transitive, or $G_u^{\Ga(u)}$ is of rank $3$ and the girth of $\Ga$ is $3$. It follows that $G/N$ is vertex-transitive on $\G_N$. If $G_u^{\Ga(u)}$ is $2$-transitive, then by part (ii), $\ov G_{\ov v}^{\Sigma(\ov v)}$ is $2$-transitive, and then by Proposition~\ref{prop-local-action}, $\Sigma$ is $(\ov G, 3)$-CH, as in part (iii). If $G_u^{\Ga(u)}$ is of rank $3$ and the girth of $\Ga$ is $3$, then by part (ii), $\ov G_{\ov v}^{\Sigma(\ov v)}$ is of rank $3$, and as $\G$ is a normal cover of $\G_N$ by part (i), the girth of $\Ga_N$ is also $3$.  Again, by Proposition~\ref{prop-local-action}, $\Sigma$ is $(\ov G, 3)$-CH, as in part (iii).\hfill\qed

\f{\bf Proof of Theorem~\ref{normal quotient}:}\ Suppose that $\Ga=(V,E)$ is a connected $(G,3)$-graph, which is not a complete multipartite graph. Then $\G$ is $G$-arc-transitive.
Let $N\lhd G$ be intransitive on $V$. If $N$ has two orbits on $V$, then since $\G$ is $G$-arc-transitive, there are no edges in each orbit of $N$. This implies that $\G$ is bipartite and hence $\G$ has girth at least $4$. By Proposition~\ref{prop-local-action}, for a vertex $u$, $G_u^{\Ga(u)}$ is $2$-transitive, and hence $\Ga$ is $(G, 2)$-arc-transitive.
If $N$ has at least three orbits on $V$, then by Lemma~\ref{quotient}~(i) and (iii), $\G_N$ is $(G/N,3)$-CH and  $\G$ is a normal cover of $\G_N$.\hfill \qed

\section{Normal $3$-CH Cayley graphs}\label{sec:normalCayley}

In this section, we investigate connected normal $3$-CH Cayley graphs. The following proposition is the main result of this section, which is useful in the study of $(G,3)$-CH graphs $\G=(V,E)$ such that $G$ is quasiprimitive on $V$ and $G$ has a minimal normal subgroup which is regular on $V$. In particular, it will be used in the proof of Theorem~\ref{th-quasiprimitive}~(1) and (3) (see Sections~\ref{sec:part1} and \ref{sec:examples}).

\begin{prop}\label{normal cayley}
Suppose that $\G=\Cay(H,S)$ is a connected Cayley graph of a group $H$ such that $\G$ is $(G,3)$-CH and of girth $3$, and $H_R\unlhd G\leq\Aut\G$. Assume that $\G$ is neither a complete graph nor a complete multipartite graph. Let $u$ be the vertex corresponding to the identity of $H$. Then $G_u\leq\Aut(H,S)$ and
there exists an integer $m$ such that $|a|=m$ for each $a\in S^2\setminus(S\cup \{u\})$. Furthermore, one of the following holds.
\begin{enumerate}
\item [{\rm (1)}]\ $S=\{s_1, s_1^{-1}, \ldots, s_r, s_r^{-1}\}$, where $r\geq 2$ and $|s_i|=3$ with $i=1,2,\ldots,r$,
$[S]\cong r\K_2$, $G_{us_1}$ is transitive on $S\setminus\{s_1,s_1^{-1}\}$ and $G_u\leq \mz_2\wr S_r$ is $2$-transitive on the set $\{\l s_1\r,\l s_2\r,\dots,\l s_r\r\}$.

\item [{\rm (2)}]\ $S$ consists of involutions, and either $H$ is an elementary abelian $2$-group, or
for any $s,s'\in S$, the product $ss'$ lies in $S$ if and only if $ss'=s's$. Moreover, either
    \begin{enumerate}
    \item [{\rm (i)}]\ $[S]$ is connected and of girth $3$, and $G_u$ is primitive on $S$ of rank $3$, or
    \item [{\rm (ii)}]\ $[S]\cong r\K_{2^n-1}$ for some positive integer $n$, and $G_u\leq S_{2^n-1}\wr S_r$.
        \end{enumerate}
      \end{enumerate}
Furthermore, if $S^H=S$ and $S$ consists of involutions, then either $H$ is a $2$-group, or $m$ is odd.
\end{prop}
\f\demo Since $H_R\unlhd G$,  the stabiliser $G_u\leq\Aut(H,S)$. Since $\G$ is a $(G,3)$-CH graph, $G_u$ acts transitively on $S=\Ga(u)$ and $\Ga_2(u)$. Thus the elements of $\Ga(u)$ are of the same order, and so are the elements of $\Ga_2(u)$. Clearly, $\G_2(u)=S^2\setminus(S\cup \{u\})$. It follows that there exists an integer $m$ such that $|a|=m$ for each $a\in S^2\setminus(S\cup \{u\})$.

Notice that $\G$ is neither a complete graph nor a complete multipartite graph and $\G$ has girth $3$. If the induced subgraph $[S]=[\Ga(u)]$ is connected, then by Theorem~\ref{local-pty}, $G_u^{\G(u)}$ is primitive of rank $3$. Since $G_u\leq\Aut(H,S)$, $S$ consists of involutions. If $[S]$ is disconnected, then by Theorem~\ref{local-pty}, the subgraph $[S]=r\K_\ell$ for some integers $r,\ell\ge2$. Since $\Ga$ is $(G,3)$-CH, $\Ga$ is $(G,2)$-geodesic transitive. For an arbitrary element $s\in S$, by \cite[Theorem~1.2]{DJLP-JACO}, we have $\l s\r\setminus\{1\}\subset S$. Let $\Delta\subset S$ be such that $s\in\Delta$ and $[\Delta]=\K_\ell$. Then $\l s\r\setminus\{1\}\subset \Delta$. Since the stabiliser $G_{us}\leq\Aut(H,S)$, $G_{us}$ fixes $s^i$ for any integer $i$. On the other hand, as $\Ga$ is $(G,3)$-CH, $G_{us}$ is transitive on $\Delta\setminus\{s\}$. We then conclude that either $|s|=2$, or $|s|=3$ and $\Delta\setminus\{s\}=\{s^{-1}\}$. As a result, we always have either $S$ consists of elements of order $3$, or $S$ consists of involutions.

Assume first that $S$ consists of elements of order $3$.
By the argument in the above paragraph, we see that $[S]=r\K_2$. Furthermore, $S=\{s_1,s_1^{-1}, s_2,s_2^{-1}, \dots, s_r,s_r^{-1}\}$ such that $\{s_i,s_i^{-1}\}$ is a block of $G_u$ acting on $\Ga(u)=S$.
By Theorem~\ref{local-pty}, $G_u^{\G(u)}=G_u^S$ is of rank 3. As $G_u\leq\Aut(H,S)$, the three orbits of $G_{us_1}$ on $S$ are $\{s_1\}, \{s_1^{-1}\}, S\setminus\{s_1,s_1^{-1}\}$. This implies that $G_u$ is 2-transitive on the set $\{\{s_1,s_1^{-1}\},\{s_2,s_2^{-1}\},\dots,\{s_r,s_r^{-1}\}\}$, and equivalently, $G_u$ is 2-transitive on $\{\l s_1\r,\l s_2\r,\dots,\l s_r\r\}$. Part (1) holds. 

Assume next that $S$ consists of involutions. Take $s\in S$. Let $m=|s|$ for some $s\in S^2\setminus (S\cup\{u\})$.
Then each element in $\Ga_2(u)=S^2\setminus (S\cup\{u\})$ has order $m$.
If $m=2$, then the product $s_1s_2$ of any elements $s_1,s_2\in S$ is of order 2, and so $s_1s_2=s_2s_1$.
It follows that $H=\lg S\rg$ is an elementary abelian $2$-group. Let $m>2$ and take $s_1, s_2\in S$.
If $s_1s_2=s_2s_1$, then the order $|s_1s_2|=2$, and so $s_1s_2\in S$, and the vertices $s_1, s_2$ are adjacent.
Conversely, if $s_1s_2\in S$ (or equivalently $s_1, s_2$ are adjacent),  then $s_1s_2$ is an involution, and hence $s_1s_2=s_2s_1$.
Consequently, $s_1s_2\in S$ if and only if $s_1s_2=s_2s_1$. This proves the first statement of part (2).
If $[S]$ is connected, then by Theorem~\ref{local-pty}, $G_u^{\G(u)}$ is primitive of rank $3$. Let $s_1,s_2\in S$ be adjacent. Then there exists $s_3\in S$ such that $s_2=s_3s_1$. As $s_3s_1\in S$, by the argument above, we have $s_3s_1=s_1s_3$. So $\{s_1, s_1s_3, s_3\}$ induces a triangle, and hence $[S]$ is of girth $3$. Part (2)(i) holds.  If $[S]$ is disconnected, then by Theorem~\ref{local-pty}, $[S]\cong r\K_\ell$ for some integers $r,\ell>1$. Then we have $S=\bigcup_{i=1}^r\Om_i$ such that $[\Om_i]\cong \K_\ell$.
The argument above shows that $\lg\Om_1\rg\cong\mz_2^n$ with $\ell=2^n-1$, as in part~(2)(ii).

Finally, we prove the last statement. Suppose that $S^H=S$ and $S$ consists of involutions of $H$.
Then $S$ is a union of full conjugacy classes of involutions of $H$, that is, $S=\bigcup_{i=1}^kC_i$, where $C_i$ is a conjugacy class of $H$.
Let $s_1,s_2\in S$, and let $m=|s_1s_2|$.
Suppose that $m$ is even.
Then $|(s_1s_2)^2|<|s_1s_2|$.
Since $S^H=S$, one has $s_2s_1s_2\in S$, and so $(s_1s_2)^2=s_1(s_2s_1s_2)\in S^2$.
Since the order $|(s_1s_2)^2|$ is less than $m=|s_1s_2|$, it follows that $s_1(s_2s_1s_2)\in S\cup\{1\}$.
It implies that $|s_1s_2|=2$ or $4$, and so the subgroup $\l s_1,s_2\r$ is of order 4 or 8.
Recall a result of Baer and Suzuki (see, for example, \cite[Theorem~2.66]{Gorenstein}) which states that if any two elements of
a conjugacy class $C$ of a finite group generate a $p$-group with $p$ prime, then all elements of $C$ are contained in a $p$-group.
We conclude that $\lg C_i\rg$ is a 2-group for all $i$.
Since $\lg C_i\rg\unlhd H$ and $\G$ is connected, one has $H=\lg S\rg=\prod_{i=1}^k\lg C_i\rg$.
Therefore, $H$ is a 2-group.
\hfill\qed

\section{A few results on $2$-transitive permutation groups}\label{sec:2-transitive}

In this section, we shall recall and prove some results on $2$-transitive permutation groups, which will be used in the proof of Theorem~\ref{th-quasiprimitive}~(1) and (2) (see Sections~\ref{sec:part1} and \ref{sec:part2}). The first is a result of Burnside (see \cite[p.202]{Burnside} or \cite[Proposition~5.2]{Cameron-2-transitive}).

\begin{prop}\label{Burnside}
A $2$-transitive group has a unique minimal normal subgroup, which is elementary abelian or simple.
\end{prop}

We will adopt the following assumption for the remainder of this section.

\medskip
\f{\bf Assumption~I.}\  $G$ is an almost simple $2$-transitive permutation group on a finite set $\Om$ and its socle is $T=\soc(G)$.

\medskip
By \cite[p.8]{Cameron-2-transitive}, $T$ is one of the groups listed in Table~\ref{table-1}.

{\begin{table}[ht]
\center
\caption{Almost simple $2$-transitive permutation groups} \label{table-1}
\begin{tabular}{llll}
\hline $T$   & point stabiliser in $T$ & degree $|\Om|$& remark \\
\hline $A_n(n\geq 5)$   & $A_{n-1}$ & $n$ &\\
$\PSL(d,q) (d\geq 2)$ & $q^{d-1}:\GL(d-1,q)/\mz_{(q-1,d)}$ & $q^d-1/q-1$ &$(d,q)\neq (2,2), (2,3)$\\

PSU$(3,q)$   & $[q^3]:\mz_{q^2-1/(3,q+1)}$ & $q^3+1$ & $q>2$\\

${}^2B_2(q)$ & $[q^2]:\mz_{q-1}$ &  $q^2+1$ & $q=2^{2a+1}>2$\\

${}^2G_2(q)$ & $[q^3]:\mz_{q-1}$ & $q^3+1$ & $q=3^{2a+1}>3$\\

$\PSp(2d,2)$ & $\GO^+(2d,2)$ & $2^{2d-1}+2^{d-1}$ & $d>2$\\

$\PSp(2d,2)$ & $\GO^-(2d,2)$ & $2^{2d-1}-2^{d-1}$ & $d>2$\\

$\PSL(2,11)$ & $A_5$ & $11$ & \\

$\PSL(2,8)$ & $D_{18}$ & $28$ & \\

$A_7$ & $\PSL(2,7)$ & $15$ &\\

$M_{11}$ & $M_{10}$ & $11$ & \\

$M_{11}$ & $\PSL(2,11)$ & $12$ & \\

$M_{12}$ & $M_{11}$ & $12$ &\\

$M_{22}$ & $\PSL(3,4)$ & $22$ & \\

$M_{23}$ & $M_{22}$ & $23$ & \\

$M_{24}$ & $M_{23}$ & $24$ & \\

$HS$ & $\PSU(3,5):\mz_2$ & $176$ & \\

$Co_3$ & $McL:\mz_2$ & $276$ & \\
\hline
\end{tabular}
\end{table}}

By inspecting Table~\ref{table-1} and \cite[p.9, Notes 2-3]{Cameron-2-transitive}, we have the following result.

\begin{prop}\label{cor-2-tran-socle}
$(1)$\ The centre of the point-stabiliser of $T$ is trivial. $(2)$\ $T$ is $2$-transitive on $\Om$ whenever $T\neq\PSL(2,8)$. $(3)$\ If $T\cong\PSL(2,8)$, then $T$ is primitive on $\Om$ and $G=\PL(2,8)\cong{}^2G_2(3)$.
\end{prop}

The following theorem gives some properties of the point stabiliser of $G$ which will be used in the proof of Theorem~1.5.

\begin{theorem}\label{lem-centralizer}
Under the notation set in Assumption~I, take $\a\in\Om$ and suppose that $G_\a$ acts $2$-transitively on a finite set $\D$ and that $T_\a$ is transitive on $\D$. Suppose further that there exists $\d\in\D$ such that the stabiliser $(T_\a)_\d$ centralises a non-identity element $t\in T$. Then $t\in (T_\a)_\d$ and $(T, T_\a, (T_\a)_\d)$ lies in {\rm Table~\ref{table-2}}.
{\footnotesize\rm\begin{table}[ht]
\caption{The list of all possible triples $(T, T_\a, (T_\a)_\d)$} \label{table-2}
\center
\begin{tabular}{llllll}
\hline No.&$T$   &$T_\a$ & $T_{\a\d}$ & $t\in$ & {\rm Remark}\\
\hline  1&$A_5$   & $A_4$ & $\mz_3$ & $T_{\a\d}$& $|t|=3$\\

2&$A_5$ & $A_4$ & $\mz_2\times\mz_2$ &$T_{\a\d}$& $|t|=2$\\

3&$A_5$ & $D_{10}$ & $\mz_5$ &$T_{\a\d}$&$|t|=5$\\

4&$\PSL(2,8)$ & $D_{18}$ & $\mz_9$ &$T_{\a\d}$& $|t|=3$ {\rm or} $9$\\

5&$\PSL(2,q)$ & $q: \mz_{k}$ & $\mz_{k}$ &$T_{\a\d}$& $k=\frac{q-1}{(2,q-1)}\neq 2$\\

 6&$\PSL(d,q)$ & $q^{d-1}: \GL(d-1,q)/\mz_k$ & $\GL(d-1,q)/\mz_k$ & $Z(T_{\a\d})$ & $d\geq 3, k=(d,q-1)$ \\

7&$\PSL(d,q)$ & $q^{d-1}: \GL(d-1,q)$ & $[q^{2d-3}]: \GL(d-2,q)$ & $q^{d-1}$ & $d\geq 3, (q-1)\ |\ d$ \\




8&$\PSU(3,q)$ & $[q^3]: \mz_{(q^2-1)/(3,q+1)}$ & $[q]:\mz_{(q^2-1)/(3,q+1)}$&$\mz_{(q+1)/(3,q+1)}$& $|t|\ |\ q+1$\\

9&$\PSU(3,3)$ & $[3^3]: \mz_{8}$ & $[3^3]: \mz_{4}$ & $Z([3^3])$ & $|t|=3$\\

10&$\PSU(3,4)$ & $[4^3]: \mz_{15}$ & $[4^3]: \mz_{5}$&$Z([4^3])$& $|t|=2$\\

11&$\PSU(3,8)$ & $[8^3]: \mz_{21}$ & $[8^3]: \mz_{3}$&$Z([8^3])$& $|t|=2$\\

12&$\PSU(3,32)$ & $[32^3]: \mz_{341}$ & $[32^3]: \mz_{11}$&$Z([32^3])$& $|t|=2$\\
\hline
\end{tabular}
\end{table}}
\end{theorem}

This theorem will be proved by the following series of lemmas. In the proof, we will adopt the following assumption.

\medskip
\f{\bf Assumption~II.}\ \begin{enumerate}
                          \item [{\rm (1)}]\ Let $K$ be the kernel of $G_\a$ acting on $\D$.
                          \item [{\rm (2)}]\ For a group $P$, we write $P^*=P\setminus\{1_P\}$, where $1_P$ is the identity element of $P$, and denote by $\Z(P)$, $P'$ and $\Phi(P)$ the centre, the derived subgroup and the Frattini subgroup of $P$, respectively.
                        \end{enumerate}


\begin{lem}\label{lem-1}
$G_\a\leq T_\a.\Out(T)$, and $|\D|-1$ divides $|T_{\a\d}/K\cap T_{\a\d}|\cdot |\Out(T)|$.
\end{lem}

\f\demo Note that $G\leq T.\Out(T)$ and $G=TG_\a$. It follows that $G_\a/T_\a\cong G/T\leq\Out(T)$. Since $T_\a$ is transitive on $\D$, one has $G_\a=T_\a G_{\a\d}$, and so $G_{\a\d}/T_{\a\d}\cong G_\a/T_\a$. As $G_\a$ is $2$-transitive on $\D$, it implies that $|\D|-1$ divides $|G_{\a\d}/K|$.
Note that $G_{\a\d}=T_{\a\d}.L$ with $L\leq\Out(T)$. So, $|G_{\a\d}/K|=\frac{|T_{\a\d}||L|}{|K|}(=\frac{|T_{\a\d}|}{|K\cap T_{\a\d}|}\cdot\frac{|L||K\cap T_{\a\d}|}{|K|})$ which divides $|T_{\a\d}/K\cap T_{\a\d}|\cdot|\Out(T)|$. Thus,
$|\D|-1$ divides $|T_{\a\d}/K\cap T_{\a\d}|\cdot|\Out(T)|$.\hfill\qed

\begin{lem}\label{lem-2}
If $T_\a\cap K\neq 1$, then $t\in T_\a$.
\end{lem}

\f\demo Clearly, $T_\a\cap K\leq T_{\a\d}$, so $T_\a\cap K$ is normal in $\lg t,T_\a\rg$ as $t$ centralises $T_{\a\d}$. By Proposition~\ref{cor-2-tran-socle}, $T$ is primitive on $\Om$, and so $T_\a$ is maximal in $T$. It follows that either $\lg t,T_\a\rg=T$ or $t\in T_\a$. If $T_\a\cap K\neq 1$, then since $T_\a\cap K\unlhd\lg t,T_\a\rg$ and $T$ is non-ableian simple, one has $t\in T_\a$.\hfill\qed

\begin{lem}\label{lem-2'}
Let $H$ be a proper abelian subgroup of $T_{\a\d}$. If $C_T(H)=H$, then $t\in H\cap Z(T_{\a\d})$.
\end{lem}

\f\demo Since $H$ is abelian, one has $H\leq C_T(H)$. If $H=C_T(H)$, then since $1\leq H\leq T_{\a\d}$ and $t\in C_{T}(T_{\a\d})$, it follows that $t\in H\leq T_{\a\d}$, and hence $t\in Z(T_{\a\d})$. So $t\in H\cap Z(T_{\a\d})$.\hfill\qed

The next lemma excludes the case where $T_\a$ is almost simple.

\begin{lem}\label{lem-3}
$T_\a$ is not almost simple.
\end{lem}

\f\demo Suppose on the contrary that $T_\a$ is almost simple. If $T_\a$ is not faithful on $\D$, then $T_\a\cap K$ is a non-trivial proper normal subgroup of $T_\a$, and so $\soc(T_\a)\leq T_\a\cap K\leq T_{\a\d}$. By inspecting Table~\ref{table-1}, we have $T_\a=\GO^\pm(2d,2),$ $M_{10}$, $\PSU(3,5):\mz_2$ or $McL:\mz_2$. By Lemma~\ref{lem-2}, we have $t\in T_\a$. Since $t$ centralises $T_{\a\d}$, it follows that $t$ centralises $\soc(T_\a)$. However, this is impossible. Thus, $T_\a$ acts faithfully on $\D$, and so $T_\a\cong T_\a K/K\unlhd G_\a/K$. Since $T_\a$ is almost simple, by Proposition~\ref{Burnside}, $G_\a/K$ is almost simple. It follows that $\soc(T_\a)\cong\soc(T_\a)K/K=\soc(G_\a/K)$. Since $G_\a/K$ is a $2$-transitive permutation group on $\D$ and $\soc(G_\a/K)\cong\soc(T_\a)\neq\PSL(2,8)$, by Proposition~\ref{cor-2-tran-socle}, $\soc(T_\a)$ is also $2$-transitive on $\D$. Inspecting the groups in Table~\ref{table-1}, we have the following possible candidates for the pair $(T, T_\a)$:
\begin{equation*}
\begin{array}{r}
    (A_n, A_{n-1})(n\geq 6), (\PSp(6,2), \GO^+(6,2)), (\PSL(2,11), A_5), (A_7, \PSL(2,7)),\\
    (M_{11}, M_{10}),(M_{11}, \PSL(2,11)), (M_{12}, M_{11}), (M_{22}, \PSL(3,4)),\\
    (M_{23}, M_{22}), (M_{24}, M_{23}), (HS, \PSU(3,5):\mz_2).
\end{array}
\end{equation*}

Let $(T, T_\a)=(A_n, A_{n-1}) (n\geq 6)$. Since $T_\a=A_{n-1}$ is a $2$-transitive permutation group on $\D$, again by checking Table~\ref{table-1}, we have either $T_{\a\d}=A_{n-2}$, or $(T, T_\a, T_{\a\d})$ is one of the following:
$$(A_6, A_5, D_{10}), (A_7, A_6, 3^2:\mz_4), (A_8, A_7, \PSL(2,7)), (A_9, A_8, 2^3:\PSL(3,2)).$$
For the former, we have $C_T(T_{\a\d})=C_{A_n}(A_{n-2})=1$ as $n\geq 6$, a contradiction. If $(T, T_\a, T_{\a\d})=(A_6, A_5, D_{10})$, then the unique Sylow $5$-subgroup of $T_{\a\d}$ is self-centralising in $T$, and by Lemma~\ref{lem-2'}, we have $t\in Z(T_{\a\d})$, but $T_{\a\d}=D_{10}$ has a trivial centre, a contradiction. If $(T, T_\a, T_{\a\d})=(A_7, A_6, 3^2: \mz_4)$, then the unique Sylow $3$-subgroup of $T_{\a\d}$ is self-centralising in $T$ and the center $T_{\a\d}=3^2: \mz_4$ is trivial. Again, this is impossible by Lemma~\ref{lem-2'}. If $(T, T_\a, T_{\a\d})=(A_8, A_7, \PSL(2,7))$ or $(A_9, A_8, 2^3: \PSL(3,2))$, then $Z(T_{\a\d})=1$. Let $P$ be a Sylow $7$-subgroup of $T_{\a\d}$. Then $C_T(P)=P$. By Lemma~\ref{lem-2'}, we obtain a contradiction.

Let $(T, T_\a)=(\PSp(6,2), \GO^+(6,2))$. Note that $\GO^+(6,2)\cong S_8\cong\PGL(4,2)$. By Table~\ref{table-1}, we have $T_{\a\d}=S_7$ or $2^3: \PSL(3,2).2$. Let $P$ be a Sylow $7$-subgroup of $T_{\a\d}$. Then $P$ is also a Sylow $7$-subgroup of $T$, and by Magma~\cite{BCP}, $C_T(P)=P$ and $Z(T_{\a\d})=1$. This is impossible by Lemma~\ref{lem-2'}.

Let $(T, T_\a)=(\PSL(2,11), A_5)$. Note that $A_5\cong\PSL(2, 5)\cong \PSL(2,4)$. By Table~\ref{table-1}, we have $T_{\a\d}=A_4$ or $D_{10}$. Then $Z(T_{\a\d})=1$. If $T_{\a\d}=A_4$, then the Sylow $2$-subgroup of $T_{\a\d}$ is self-centralising in $T$, and if $T_{\a\d}=D_{10}$, then the Sylow $5$-subgroup of $T_{\a\d}$ is self-centralising in $T$. By Lemma~\ref{lem-2'}, we see that this case cannot happen.

Let $(T, T_\a)=(A_7, \PSL(2,7))$. Note that $\PSL(2,7)\cong\PSL(3,2)$. By Table~\ref{table-1}, we have $T_{\a\d}=\mz_7:\mz_3$ or $S_4$. Then $Z(T_{\a\d})=1$. If $T_{\a\d}=\mz_7:\mz_3$, then the Sylow $7$-subgroup of $T_{\a\d}$ is self-centralising in $T$, which is impossible by Lemma~\ref{lem-2'}. If $T_{\a\d}=S_4$, then by Magma~\cite{BCP}, we obtain that $C_T(T_{\a\d})=1$, a contradiction.

Let $(T, T_\a)=(M_{11}, M_{10})$. Note that $\soc(M_{10})=A_6$. By Table~\ref{table-1}, we have $T_{\a\d}=A_5.2$ or $\mz_3^2: D_8$. However, by Magma~\cite{BCP}, $M_{10}$ has no subgroups of order $120$, and if $T_{\a\d}=\mz_3^2: D_8$, then $C_T(T_{\a\d})=1$, a contradiction.

Let $(T, T_\a)=(M_{11}, \PSL(2,11))$. By Table~\ref{table-1}, we have $T_{\a\d}=\mz_{11}: \mz_5$.  Then $Z(T_{\a\d})=1$ and the Sylow $11$-subgroup of $T_{\a\d}$ is self-centralising in $T$, which is impossible by Lemma~\ref{lem-2'}.

Let $(T, T_\a)=(M_{12}, M_{11})$. By Table~\ref{table-1}, we have $T_{\a\d}=M_{10}$ or $\PSL(2,11)$.  However, by Magma~\cite{BCP}, we obtain that $C_T(T_{\a\d})=1$, a contradiction.

Let $(T, T_\a)=(M_{22}, \PSL(3,4))$. By Table~\ref{table-1}, we have $T_{\a\d}=4^2:\PGL(2,4)$. However, by Magma~\cite{BCP}, we obtain that $C_T(T_{\a\d})=1$, a contradiction.

Let $(T, T_\a)=(M_{23}, M_{22})$. By Table~\ref{table-1}, we have $T_{\a\d}=\PSL(3,4)$. However, by Magma~\cite{BCP}, we obtain that $C_T(T_{\a\d})=1$, a contradiction.

Let $(T, T_\a)=(M_{24}, M_{23})$. By Table~\ref{table-1}, we have $T_{\a\d}=M_{22}$. However, by Magma~\cite{BCP}, we obtain that $C_T(T_{\a\d})=1$, a contradiction.

Let $(T, T_\a)=(HS, \PSU(3,5):\mz_2)$. By Table~\ref{table-1}, we have $T_{\a\d}=[5^3]:\mz_{8}.2$. However, by Magma~\cite{BCP}, we obtain that $C_T(T_{\a\d})=1$, a contradiction. \hfill\qed

Now, we consider the case where $T_\a$ is soluble and $T\neq\PSL(d,q) (d\geq 3)$.

\begin{lem}\label{lem-4}
If $T_\a$ is solvable and $T\neq\PSL(d,q) (d\geq 3)$, then lines~{\rm 1--5}, and lines~{\rm 8--12} of {\rm Table~\ref{table-2}} hold.
\end{lem}

\f\demo By Lemma~\ref{lem-1}, $G_\a$ is solvable, and so $\soc(G_\a/K)\leq T_\a K/K$ and is elementary abelian. By Table~\ref{table-1}, we have the following possible candidates for $(T, T_\a)$:
\begin{equation*}
 \begin{array}{l}
 (A_5, A_4), (\PSL(2,q), [q]:\mz_{q-1/(q-1,2)}), (\PSU(3, q), [q^3]:\mz_{q^2-1/(3,q+1)}), \\ ({}^2B_2(q), [q^2]:\mz_{q-1}),
  ({}^2G_2(q), [q^3]:\mz_{q-1}), (\PSL(2,8), D_{18}).
 \end{array}
\end{equation*}

Let $(T, T_\a)=(A_5, A_4)$. If $T_{\a}\cap K=1$, then the unique minimal normal subgroup of $T_\a$ is regular on $\D$, and hence
$T_{\a\d}\cong\mz_3$. As $T_{\a\d}$ is self-centralising in $T$ and $t\in C_T(T_{\a\d})$, one has $t\in T_{\a\d}$, as in line~1 of Table~\ref{table-2}. If $T_{\a}\cap K\neq 1$, then $T_{\a\d}=T_{\a}\cap K\cong\mz_2^2$ and $T_\a K/K\cong\mz_3$ is regular on $\D$. Similarly, since $T_{\a\d}$ is self-centralising in $T$  and $t\in C_T(T_{\a\d})$, one has $t\in T_{\a\d}\cong\mz_2^2$, as in line~2 of Table~\ref{table-2}.

Let $(T, T_\a)=(\PSL(2,8), D_{18})$. Since $T_\a\unlhd G_\a$, the unique cyclic subgroup of $T_\a$ of order $9$ is normal in $G_\a$. Considering the $2$-transitivity of $G_\a$ on $\D$ shows that $T_{\a}\cap K\neq 1$. By Lemma~\ref{lem-2}, we have $t\in T_\a$, and since $T_\a\cong D_{18}$ and $t$ centralises $T_{\a\d}$, one has $T_{\a\d}=T_{\a}\cap K\cong\mz_9$ and $t\in T_{\a\d}$, as in line~4 of Table~\ref{table-2}.

Let $(T, T_\a)=(\PSL(2,q), \mz_p^r: \mz_{(q-1)/d})$, where $q=p^r$, $d=(q-1,2)$ and $q\neq 2,3$. Note that $T_\a=\mz_p^r: \mz_{(q-1)/d}$ is a Frobenius group and $\mz_p^r$ is minimal normal in $T_\a$.  Assume first that $T_{\a}\cap K=1$, then $T_{\a\d}=\mz_{(q-1)/d}$ and $C_T(T_{\a\d})=T_{\a\d}$. So, $t\in T_{\a\d}$. If $(q-1)/d\neq 2$, then we obtain line~5 of Table~\ref{table-2}. If $(q-1)/d=2$, then $q$ is odd and so $d=(q-1,2)=2$. It implies that $q-1=2d=4$ and so $q=5$. Then $T\cong A_5$, as line~4 of Table~\ref{table-2}. Assume now that $T_{\a}\cap K\neq1$. Then $\mz_p^r\leq K$. It follows that $T_\a K/K$ is cyclic and hence regular on $\D$. Thus, $\soc(G_\a/K)=T_\a K/K\cong\mz_\ell$ for some prime $\ell$, and $T_\a\cap K=T_{\a\d}$. By Lemma~\ref{lem-2}, we have $t\in T_\a$. Since $t$ centralises $\mz_p^r$ and since $T_\a$ is a Frobenius group, it follows that $T_{\a\d}=T_\a\cap K=\mz_p^r$, and so $t\in\mz_p^r$ and $\mz_{(q-1)/d}\cong\mz_\ell$. 
By Lemma~\ref{lem-1}, we have $|\D|-1=\ell-1$ divides $|\Out(T)|$. First, let $q=p^r$ be odd. Then $\ell=(p^r-1)/2$ and $|\Out(T)|=2r$. If $r>1$, then $(p^{r-1}+\cdots+p)\ |\ 2r$, which is impossible. If $r=1$, then $(p-1)/2$ divides $2$ and then $p=5$. So $T=\PSL(2,5)\cong A_5$ and line~3 of Table~\ref{table-2} happens.
Now let $q=2^r$. Then $\ell=2^r-1$, and then $|\Out(T)|=r$ and so $(2^r-2)\ |\ r$. It follows that $r=2$ and $T=\PSL(2,4)\cong A_5$. Again, we have line~2 of Table~\ref{table-2}.

Let $(T, T_\a)=(\PSU(3, q), [q^3]:\mz_{(q^2-1)/d})$, where $d=(q+1,3)$ and $q=p^r$ for some prime $p$ and integer $r$. Set $P=[q^3]$ and $Q=\mz_{(q^2-1)/d}$ so that $T_\a=P: Q$. Let $Z=\Z(P)$. Then $T_\a$ has the following properties:
(1)\ $Z=P'=\Phi(P)$ has order $q$ (\cite[Lemma~1.13(ii)]{ONan}); (2)\ $Q$ acts, by conugation, irreducibly on $P/Z$(\cite[Lemma~1.12(ii)]{ONan}); (3)\ $Q$ acts, by conugation, faithfully and semiregularly on $P\setminus Z$. Furthermore, $C_Q(Z)=C_Q(z)\cong \mz_{(q+1)/d}$ for any $1\neq z\in Z$, and $Q/C_Q(Z)$ acts regularly by conugation on $Z^*$ (\cite[Lemma 1.14]{ONan}).

Notice that $Z$ is characteristic in $T_\a$ and so normal in $G_\a$ because $T_\a\unlhd G_\a$. If $Z\nleq K$, then $ZK/K$ is regular on $\D$, and so $Z$ is transitive on $\D$. Consequently, we have $P=ZP_{\d}=P_\d$ since $Z=\Phi(P)$, a contradiction.  Thus, $Z\leq K$. By Lemma~\ref{lem-2}, we have $t\in T_{\a}$. Since $t$ centralises $T_{\a\d}$, it also centralises $Z$ because $Z\leq K\cap T_{\a\d}$.

If $P\nleq K$, then since $P'=Z\leq K$, $PK/K$ is abelian and so regular on $\D$. Since $Q$ acts irreducibly on $P/Z$, one has $P\cap K=Z$ and so $q^2=|P/Z|=|PK/K|=|\D|$. Then $T_\a=PT_{\a\d}$ and $|T_{\a\d}|=|T_\a|/q^2=q|Q|$. Note that $P,Z\unlhd T_\a$ and all subgroups of $T_\a$ of order $(q^2-1)/d$ are conjugate. We may assume that $Q\leq T_{\a\d}$. Then $T_{\a\d}=Z: Q$, and then $t$ centralises $Q$. Since $C_{T_\a}(Q)=Q$, one has $t\in Q$. Since $t$ centralises $Z$, one has $t\in C_Q(Z)\cong \mz_{(q+1)/d}$. This is line~8 of Table~\ref{table-2}.

If $P\leq K$, then $P\leq T_{\a\d}$, and so $t\in C_{T_\a}(P)$.  If $C_{T_\a}(P)\nleq P$, then $|C_{T_\a}(P)|$ has a prime divisor $\ell$ such that $\ell\neq p$. Take an element $g\in C_{T_\a}(P)$ of order $\ell$. Clearly, $\ell\mid |Q|$, so $g\in Q^x$ for some $x\in T_\a$. Without loss of generality, we may assume that $g\in Q$. However, since $Q$ acts, by conugation, faithfully and semiregularly on $P\setminus Z$ and irreducibly on $P/Z$, it follows that $g=1$, a contradiction. Thus, $C_{T_\a}(P)\leq P$. Then $t\in P$, and then $t\in\Z(T_{\a\d})$. Furthermore, $T_\a K/K=QPK/K=QK/K$. Since $Q$ is cyclic, one has $QK/K=\soc(G_\a/K)\cong\mz_\ell$ for some prime $\ell$. It follows that $T_\a\cap K=T_{\a\d}$ which is maximal in $T_\a$. Since the centre of $T_\a$ is trivial, one has $T_{\a\d}=C_{T_\a}(t)$, and so $C_Q(t)=Q\cap T_{\a\d}=Q\cap K$.
By the above property (3), we have $C_Q(t)=C_Q(Z)$. Then $\mz_{q-1}\cong Q/C_Q(Z)=Q/C_Q(t)\cong QK/K\cong\mz_\ell$. It follows that $\ell=q-1$, and since $\ell$ is prime, one has $q=3$ or $2^r$. If $q=3$, then we have line~9 of Table~\ref{table-2}. If $q=2^r$, then by Lemma~\ref{lem-1}, we have $\ell-1\ |\ |\Out(\PSU(3,q))|$, namely, $2^r-2\ |\ 2dr$. This implies that $T=\PSU(3,q)$ with $q=4,8,$ or $32$, as in lines~10-$12$ of Table~\ref{table-2}.

Let $(T, T_\a)=({}^2B_2(q), [q^2]:\mz_{q-1})$, where $q=2^{2n+1}>2$. For convenience, we let $T_\a=P: Q$, where $P=[q^2]$ and $Q\cong \mz_{q-1}$. Then $T_\a$ has the following properties (see \cite{Suzuki-simple} or \cite[Lemma~3.2]{Fang-Praeger-Suzuki}): (1)\ $T_\a$ is a Frobenius group; (2)\ $Q$ acts transitively and regularly by conjugation on both $\Z(P)^*$ and $(P/\Z(P))^*$; (3)\ $\Z(P)$ and $P/\Z(P)$ both are elementary abelian group of order $q$, and in particular, $\Phi(P)\leq \Z(P)$.

Clearly, $\Z(P)$ is characteristic in $T_\a$. Since $T_\a\unlhd G_\a$, one has $\Z(P)\unlhd G_\a$.
If $\Z(P)\nleq K$, then $\Z(P)K/K$ is regular on $\D$, and so $\Z(P)$ is transitive on $\D$. Consequently, we have $P=\Z(P)P_{\d}$. Since $\Phi(P)\leq\Z(P)$, $P_\d\cong P/\Z(P)$ is elementary abelian, forcing that $P$ is elementary abelian, a contradiction. Thus, $\Z(P)\leq K$. By Lemma~\ref{lem-2}, we have $t\in T_{\a}$. Then $t$ centralises $\Z(P)$ because $\Z(P)\leq K\cap T_{\a\d}$. Since $Q$ acts transitively and regularly by conjugation on $\Z(P)^*$, one has $t\in P$. Since $T_\a$ is a Frobenius group, one has $T_{\a\d}\leq C_{T_\a}(t)\leq P$. If $P\nleq K$, then $P$ is transitive on $\D$ and so $T_\a=PT_{\a\d}=P$, a contradiction. Thus, $P\leq K$ and so $P\leq T_\a\cap K\leq T_{\a\d}$. As we have already shown that $T_{\a\d}\leq P$, it follows that $P=T_\a\cap K=T_{\a\d}$. Consequently, $T_\a K/K=QPK/K=QK/K\cong Q$, which is cyclic and regular on $\D$. So, $QK/K$ is the socle of $G_\a/K$, and hence $Q\cong QK/K\cong\mz_\ell$ for some prime $\ell$. It follows that $q-1=2^{2n+1}-1=\ell$. By Lemma~\ref{lem-1}, we have $\ell-1\ |\ |\Out({}^2B_2(q))|$, namely, $2^{2n+1}-2\ |\ 2n+1$, a contradiction.

At last, let $(T, T_\a)=({}^2G_2(q), [q^3]:\mz_{q-1})$, where $q=3^{2n+1}>3$. Set $T_\a=P: Q$, where $P=[q^3]$ and $Q\cong \mz_{q-1}$. Then $T_\a$ has the following properties (see, for example, \cite[Lemma~2.1]{Fang-Praeger-Ree}): (1)\ $\Z(P)<\Phi(P)=P'<P$, with $P'$ elementary abelian of order $q^2$ and $|\Z(P)|=q$; (2)\ $Q$ acts regularly by conjugation on $\Z(P)^*$; (3)\ $Q$ is self-centralising in $T_\a$.

Clearly, $P'$ is characteristic in $T_\a$. Since $T_\a\unlhd G_\a$, one has $P'\unlhd G_\a$.
If $P'\nleq K$, then $P'K/K$ is regular on $\D$, and so $P'$ is transitive on $\D$. Consequently, we have $P=P'P_{\d}=P_\d$, a contradiction. Thus, $P'\leq K$. By Lemma~\ref{lem-2}, we have $t\in T_{\a}$. Then $t$ centralises $P'$ because $P'\leq K\cap T_{\a\d}$. Since $Q$ acts regularly by conjugation on $\Z(P)^*$, one has $t\in P$. If $P\leq K$, then $t$ centralises $P$ and so $t\in\Z(P)$. Since $Q$ is cyclic, $T_\a K/K=QPK/K=QK/K$ is the socle of $G_\a/K$. Then $T_\a K/K=QK/K\cong\mz_\ell$ for some prime $\ell$. Since $t\in\Z(P)$, one has $Q\cap K=1$ by the above property (2), and so $Q\cong QK/K\cong\mz_\ell$. However, $q-1=3^{2n+1}-1$ is not a prime, a contradiction. Thus, $P\nleq K$, and so $PK/K$ is the socle of $G_\a/K$. In particular, we may assume that $Q\leq T_{\a\d}$. So, $t$ centralises $Q$. This is impossible because $Q$ is self-centralising in $T_\a$ and $t\in P$.\hfill\qed

At last, we deal with the case where $T=\PSL(d,q) (d\geq 3)$.

\begin{lem}\label{lem-5}
If $T=\PSL(d,q) (d\geq 3)$, then lines~{\rm 6--7} of {\rm Table~\ref{table-2}} hold.
\end{lem}

\f\demo Let $F_q$ be the field of order $q$, and let $\Om$ be the set of $1$-dimensional subspaces of $F_q^d$. Then $\SL(d,q)$ acts on $\Om$ with kernel $\mz:=\Z(\SL(d,q))$. Then $T=\SL(d,q)/\mz$, and $T$ is a $2$-transitive permutation group on $\Om$. Let $\a$ be the subspace $\l(1,0,\ldots, 0)\r$.

Set
$$
\begin{array}{l}
P=\{\left(
      \begin{array}{cc}
        1 & 0 \\
        u & I_{d-1} \\
      \end{array}
    \right)\ |\ u\ {\rm is\ a}\ (d-1)\times 1\ {\rm matrix\ over}\ F_q\},
    \end{array}
$$
where $I_{d-1}$ is the $(d-1)\times (d-1)$ identity matrix over $F_q$.
Then $P$ is an elementary abelian group of order $q^{d-1}$, and $P\cap \mz=1$. It is easy to see that $P$ fixes $\a$. So $P\cong P\mz/\mz\leq T_\a$.

Set
$$
\begin{array}{l}
H=\{\left(
      \begin{array}{cc}
        a & 0 \\
        0 & B \\
      \end{array}
    \right)\ |\ a\in F_q, B\in \GL(d-1,q), a{\rm det}(B)=1\},\\
L=\{\left(
      \begin{array}{cc}
        1 & 0 \\
        0 & B \\
      \end{array}
    \right)\ |\ B\in \SL(d-1,q), {\rm det}(B)=1\}.
    \end{array}
$$
Then $H\cong\GL({d-1},q)$, $L\cong\SL(d-1,q)$, and $L\leq H$. It is also easy to see that $H$ fixes $\a$ and $\mz\leq H$. So $H/\mz\leq T_\a$. Furthermore, we have $|H\cap P|=1$, and so $|PH/\mz|=|T_\a|$. Consequently, we have $PH/\mz=T_\a$.


For any $(d-1)\times 1$ matrix $u$, write $g_u=\left(\begin{array}{cc}
        1 & 0 \\
        u & I_{n-1} \\
      \end{array}\right)$.
Set $V=\{u\ |\ g_u\in P\}$. Then $V$ is an $(n-1)$-dimensional vector space over $F_q$. Define the action of $H$ on $V$ by
$$u\mapsto aB^{-1}u, \forall u\in V, h=\left(
      \begin{array}{cc}
        a & 0 \\
        0 & B \\
      \end{array}
    \right)\in H.$$
It is easy to see that this action is permutation equivalent to the conjugate action of $H$ on $P$. So, we may identify $P$ with $V$. Then $L$ acts irreducibly on $V$.
If $\left(\begin{array}{cc}
        a & 0 \\
        0 & B \\
      \end{array}
    \right)\in C_H(P)$, then we have $aB^{-1}=I_{d-1}$, and hence $B=aI_{d-1}$. This implies that $\left(\begin{array}{cc}
        a & 0 \\
        0 & B \\
      \end{array}
    \right)\in\mz$.

So we have the following claim.

\medskip
\f{\bf Claim~1.}\ $C_H(P)\leq\mz$ and $L$ acts by conjugation irreducibly and faithfully on $P$.
\medskip

Next we prove several other claims.

\medskip
\f{\bf Claim~2.}\ For arbitrary $h\in H$ and $g, g'\in P$, $g^h=g'$ if and only if $g^h\mz=g'\mz$.
\medskip

Clearly, $g^h\mz=g'\mz$ whenever $g^h=g'$. Suppose that $g^h\mz=g'\mz$. Then $(g')^{-1}g^h\in\mz\cap P=1$, and so $g^h=g'$. This completes the proof of Claim~2.

\medskip
\f{\bf Claim~3.}\ $H/\mz$ acts by conjugation irreducibly on $P\mz/\mz$.
\medskip

As $\mz\cong\mz_{(q-1,d)}$ and $\Z(\SL(d-1,q))\cong\mz_{(q-1,d-1)}$, one has  $L\cap \mz=1$, and hence
$L\cong L\mz/\mz\leq H/\mz$. By Claim~1, $L$ acts by conjugation irreducibly on $P$. If $P\mz/\mz$ has a subgroup, say $B\mz/mz$, which is normalised by $L\mz/\mz$, then for any $a\in L$, we have $B^a\mz/\mz=B\mz/\mz$ and then $B^a=B$ by Claim~2. The arbitrariness of $a$ implies that $B$ is normalised by $L$. As $L$ acts by conjugation irreducibly on $P$, we have either $B=1$ or $B=P$. It follows that $L\mz/\mz$ acts by conjugation irreducibly on $P\mz/\mz$.

\medskip
\f{\bf Claim~4.}\ $C_{T}(H/\mz)\leq H/\mz$.
\medskip

Clearly, $C_{T}(H/\mz)\leq C_T(L\mz/\mz)$ as $L\leq H$. Take $g\in C_{T}(H/\mz)$. Assume that $g=\left(
      \begin{array}{cc}
        c & u \\
        v & D \\
      \end{array}
    \right)\mz$, where $c\in F_q$ and $D$ is a $(d-1)\times (d-1)$ matrix. Then for any $B\in \SL(d-1,q)$, we have
    \[\left(
      \begin{array}{cc}
        c & u \\
        v & D \\
      \end{array}
    \right)\left(
      \begin{array}{cc}
        1 & 0 \\
        0 & B \\
      \end{array}
    \right)\mz=\left(
      \begin{array}{cc}
        1 & 0 \\
        0 & B \\
      \end{array}
    \right)\left(
      \begin{array}{cc}
        c & u \\
        v & D \\
      \end{array}
    \right)\mz.\]
It follows that \[\left(
      \begin{array}{cc}
        c & uB \\
        v & DB \\
      \end{array}
    \right)=k\left(
      \begin{array}{cc}
        c & u \\
        Bv & BD \\
      \end{array}
    \right),\]
where $k$ is a non-zero element of $F_q$. So we have $c=kc$, $DB=kBD$, $kBv=v$ and $uB=ku$. If $v\neq 0$, then $v$ is an eigenvector of $B$. By the arbitrariness of $B$, $v$ is also an eigenvector of $I_{ij}$, where $I_{ij}\in\SL(d-1,q)$ is obtained from $I_{d-1}$ by adding the $i$th row to the $j$th row. Note that all eigenvalues of $I_{ij}$ are $1$. So for $j=1,2,\dots, d-1$, we have $I_{dj}v=v$ and hence the $j$th entry of $v$ is $0$. Also, $I_{1d}v=v$ and hence the $d$th entry is $0$. This implies that $v=0$, a contradiction.  If $u\neq 0$, then $uB=ku$ implies that $B^Tu^T=ku^T$, and hence $u^T$ is an eigenvector of $B^T$. With a similar argument above, we would obtain $u=0$, a contradiction. Thus, we have $v=0$ and $u=0$. So $g\in H/\mz$. This completes the proof of Claim~4.\medskip

Now we are ready to finish the proof of this lemma.
Recall that $G_\a$ acts $2$-transitively on a finite set $\D$ with kernel $K$, and that $T_\a$ is transitive on $\D$.
Suppose that $P\mz/\mz\nleq K$. Note that $P\mz/\mz$ is characteristic in $T_\a$. Since $T_\a\unlhd G_\a$, $P\mz/\mz$ is normal in $G_\a$. Then $P\mz/\mz$ is transitive on $\D$ as  $G_\a$ acts $2$-transitively on $\D$. By Claim~3, $L\mz/\mz$ acts by conjugation irreducibly on $P\mz/\mz$. Then $(P\mz/\mz)\cap K$ is trivial, and so $P\mz/\mz$ is regular on $\D$ as $P\mz/\mz$ is elementary abelian. So, we may identify $\D$ with $P\mz/\mz$ and take $\d$ to be the identity element of $P\mz/\mz$. Then $T_\a=P\mz/\mz: T_{\a\d}$ with $T_{\a\d}=H/\mz$. This implies that $t$ centralises $H/\mz$.  By Claim~4, we have $C_{T}(H/\mz)\leq H/\mz$. It follows that $t\in\mz(H/\mz)$. We get line 6 of Table~\ref{table-2}.

Next, suppose that $P\mz/\mz\leq K$. By Lemma~\ref{lem-2}, $t\in T_{\a}$ and $t$ centralises $P\mz/\mz$. By Claim~1, we have $C_H(P)\mz/\mz=\mz$. So, $t\in P\mz/\mz$. Assume that $t=g_u\mz$ for some $0\neq u\in V$.

Suppose that $\Z(H/\mz)\cap K\neq 1$. Take $1\neq h\in \Z(H/\mz)\cap K$. Since $\mz\leq\Z(H)$, one has $\Z(H/\mz)=\Z(H)/\mz$. Then $h=h'\mz$, where $h'\in\Z(H)\setminus\mz$. So $h'=\left(
      \begin{array}{cc}
        a & 0 \\
        0 & B \\
      \end{array}
    \right)$, where $B$ is a scalar of $\GL(d-1,q)$. Then $(h')^{-1}g_uh'\mz=\left(
      \begin{array}{cc}
        a & 0 \\
        aB^{-1}u & I_{d-1} \\
      \end{array}
    \right)\mz=\left(\begin{array}{cc}
        a & 0 \\
        u & I_{d-1} \\
      \end{array}
    \right)\mz=g_u\mz$. By Claim~2, we have $aB^{-1}u=u$. Since $u\neq 0$ and $B$ is a scalar of $\GL(d-1,q)$, one has $B=aI_{d-1}$. So $h'\in\mz$, a contradiction. Thus, we have $\Z(H/\mz)\cap K=1$.

If $\Z(H/\mz)\neq 1$, then $\Z(H/\mz)\cong \Z(H/\mz)K/K$ acts regularly on $\D$. Consequently, $H/\mz=\Z(H/\mz)\times (H/\mz)_\d$, and so $(H/\mz)_\d\leq K$. Recall that $L\mz/\mz\cong\SL(d-1,q)$. As $(H/\mz)/(H/\mz)_\d$ is abelian, the derived subgroup of $H/\mz$ is contained in $(H/\mz)_\d$. Since $H\cong\GL(d-1,q)$, one has $\SL(d-1,q)\cong L\cong L\mz/\mz\leq (H/\mz)_\d$. By Claim~1, $L\mz/\mz$ is transitive on the 1-dimensional subspaces of $V$. However, $(H/\mz)_\d$ centralises $t=g_u\mz$, and so for any $x\mz\in (H/\mz)_\d$ with $x\in H$, we have $g_u^x\mz=g_u\mz$. As $g_u\in P$, Claim~2 implies that $g_u^x=g_u$. By the argument above Claim~1, we see that $x$ fixes $u$, and hence $(H/\mz)_\d$ fixes the subspace $\lg u\rg$. A contradiction occurs.

Thus, $\Z(H/\mz)=1$. Again, since $\mz\leq\Z(H)$ and $\Z(H)\cong\mz_{q-1}$, one has $1=\Z(H/\mz)=\Z(H)/\mz\cong \mz_{q-1}/\mz_{(q-1,d)}$, implying $(q-1)\ |\ d$. Thus, $\mz=\Z(H)$ and $\Z(L)=\Z(\SL(d-1,q))=1$. So $H/\mz\cong\PGL(d-1,q)$ and $\SL(d-1,q)=\PSL(d-1, q)$.  Since $|\SL(d-1,q)|=|\PGL(d-1,q)|$, one has $\PSL(d-1,q)=\SL(d-1,q)\cong L\cong L\mz/\mz=H/\mz\cong\PGL(d-1,q)$. Thus, $K=P\mz/\mz$. If $\PSL(d-1,q)$ is not simple, then since $(q-1)\mid d$ and $d\geq 3$, one has $d-1=2$ and $q=2$, and then $T_\a K/K\cong H/\mz\cong\PGL(2,2)\cong S_3$. It follows that $(H/\mz)_\d\cong\mz_2$ as line~7 of Table~\ref{table-2}. Assume now that $\PSL(d-1, q)$ is simple. Then $T_\a K/K\cong H/\mz$ is the socle of $G_\a/K$. Since $G_\a/K$ is $2$-transitive on $\D$, by Proposition~\ref{cor-2-tran-socle}, $H/\mz$ is also $2$-transitive on $\D$ and $T_{\a\d}=(P\mz/\mz): (H/\mz)_\d$. By inspecting Table~\ref{table-1} (see also \cite[Main theorem]{Curtis-2-tran-lie}) and using the fact that $(q-1)\mid d$, we have either $(H/\mz)_{\d}=q^{d-2}: \GL(d-2,q)$, or $(H/\mz, (H/\mz)_\d)$ is one of the followings:
\[(\PSL(2, 4), A_4), (\PSL(3, 2), \mz_7: \mz_3), (\PSL(4, 2), A_7).\]
If the latter happens, then by Magam~\cite{BCP} we obtain that $C_T(T_{\a\d})$ is trivial, contrary to $t\in C_T(T_{\a\d})$.
Thus, line~7 of Table~\ref{table-2} happens.\hfill\qed

\medskip
\f{\bf Proof of Theorem~\ref{lem-centralizer}}\ By Lemma~\ref{lem-3}, $T_\a$ is not almost simple. If $T_\a$ is solvable, then by Lemma~\ref{lem-4}, lines 1--5 or 8--12 of Table~\ref{table-2} happen. If $T_\a$ is neither almost simple nor solvable, then $T=\PSL(d,q) (d\geq 3)$, and by Lemma~\ref{lem-5}, lines~6--7 of Table~\ref{table-2} hold. \hfill\qed



\section{There exist no quasiprimitive $3$-CH graphs of type HS, HC and CD}\label{sec:part1}

In this section, we prove the following theorem, which is just part (1) of Theorem~\ref{th-quasiprimitive}.

\begin{theorem}\label{th:hs-hc-cd}
Let $\Ga=(V,E)$ be a connected $(G,3)$-CH graph of girth $3$, which is neither complete nor complete multipartite. If $G$ is quasiprimitive on $V$, then $G$ is not of type {\rm HS, HC} or {\rm CD}.
\end{theorem}

\f\demo Assume that $G$ is quasiprimitive on $V$. We shall finish the proof by the following two steps.


\medskip
\f{\bf Step~1.}\ $G$ is not of type {\rm HS} or {\rm HC}.

Suppose on the contrary that $G$ is quasiprimitive on $V$ of type {\rm HS} or {\rm HC}.
Then $G$ has a non-abelian minimal normal subgroup, say $H$, which is regular on $V\G$, and
the stabiliser $G_u$ with $u$ being the identity of $H$ is such that $\Inn(H)\leq G_u\leq \Aut(H)$.
The subgroup $H=T^k$, where $k$ is an integer and $T$ is a non-abelian simple group.
We may view $\G$ as a Cayley graph $\Cay(H, S)$ on the group $H$.
Then the neighbourhood $S:=\Ga(u)$ is a subset of $H$, and $\Inn(H)\lhd G_u\le\Aut(H,S)$.
In particular, $S$ is a union of full conjugacy classes of elements of $H$. By Proposition~\ref{normal cayley}, all elements of $S$ have order $n$, where $n=3$ or $2$.

Suppose first that all elements of $S$ have order $3$.
By Proposition~\ref{normal cayley}~(1), $[S]\cong r\K_2$, $S=\{s_1,s_1^{-1},\dots,s_r,s_r^{-1}\}$, where $r\geq 2$ and
$G_u$ is 2-transitive on the set $\{\l s_1\r,\l s_2\r,\dots,\l s_r\r\}$. Since each $s_i$ has order $3$, $G_u$ is also 2-transitive on $\calB=\{\Delta_i\mid \D_i=\{s_i,s_i^{-1}\}, i=1,2,\dots,r\}$.
Let $K$ be the kernel of $G_u$ acting on $\calB$.
Then $K$ is a 2-group, and $G_u/K$ is a 2-transitive permutation group on $\calB$. As $\Inn(H)\lhd G_u\leq\Aut(H)$, one has $\Inn(H)K/K\lhd G_u/K$, and as $K$ is a $2$-group and $\Inn(H)\cong T^k$ with $T$ a non-abellian simple group, it follows that $K\cap\Inn(H)=1$ and hence $\Inn(H)\cong\Inn(H)K/K\leq\soc(G_u/K)$. Since $G_u/K$ is a 2-transitive permutation group on $\calB$, by Proposition~\ref{Burnside} we see that $k=1$ and $\soc(G_u/K)\cong H=T$. Thus, $G$ is of type HS. Noticing $\Inn(H)\lhd G_u\leq\Aut(H)$, $G_u$ is almost simple with socle $\Inn(H)\cong T$. This implies that $K=1$ as $K$ is a normal $2$-subgroup of $G_u$.
In particular, $G_u$ is a 2-transitive permutation group on $\calB$.
Recall that $\soc(G_u)=\Inn(T)$. By Proposition~\ref{cor-2-tran-socle}~(2)-(3), $\Inn(T)$ is primitive on $\calB$ and so $\Inn(T)_{\D_1}$ is maximal in $\Inn(T)$. If $\Inn(T)_{\D_1}$ is intransitive on $\D_1=\{s_1, s_1^{-1}\}$, then $\Inn(T)_{\D_1}$ centralises $s_1$. The maximality of $\Inn(T)_{\D_1}$ implies that $\s(s_1)\in \Inn(T)_{\D_1}$, where $\s(s_1)$ is the inner automorphism of $T$ induced by $s_1$. So the centre of $\Inn(T)_{\D_1}$ is nontrivial, which is impossible by Proposition~\ref{cor-2-tran-socle}~(1). Thus, $\Inn(T)_{\D_1}$ is transitive on $\D_1=\{s_1, s_1^{-1}\}$.
Now apply Theorem~\ref{lem-centralizer} to the 2-transitive permutation group $G_u^{\calB}$ on $\calB$, and the transitive action of $\Inn(T)_{\D_1}$ on $\D_1=\{s_1, s_1^{-1}\}$. Note that $\Inn(T)_{\D_1 s_1}$ centralises $\s(s_1)$. Then $\Inn(T), \Inn(T)_{\D_1}, \Inn(T)_{\D_1 s_1}$ and $\s(s_1)$ satisfy Table~\ref{table-2}. Since $\Inn(T)_{\D_1 s_1}$ has index $2$ in $\Inn(T)_{\D_1}$ and $|s_1|=3$, the lines~4 or 9 of Table~\ref{table-2} happen. Then $H=T=\PSL(2,8)$ or $\PSU(3,3)$, $|\calB|=28$, $|S|=56$ and $\G=\Cay(H,S)$. By Magma~\cite{BCP}, $S$ is a conjugacy class of $H$ of elements of order $3$. However, by Magma~\cite{BCP}, we see that $\Aut(\G)_u^S$ is of rank $4$, which is impossible by Proposition~\ref{prop-local-action}.

Suppose now that all elements of $S$ have order $2$.
By Proposition~\ref{normal cayley}~(2), for any $s_1, s_2\in S$, either $s_1s_2\in S\cup\{u\}$ or $o(s_1s_2)=m$ for a fixed odd integer $m$.
Recall that $S$ is a union of full conjugacy classes of involutions of $H$. Let $D\subset S$ be a conjagacy class of involutions of $H$. Let $L=\lg D\rg$. Then $L\unlhd H$, and so $L=T^\ell$ for some $1\leq \ell\leq k$. By \cite[Theorem]{Aschbacher}, if $B$ is a group which has no nontrivial solvable normal subgroups and $B'=B''$, and if $B$ is generated by a conjugacy class $D$ of transpositions such that the product of any two distinct members of $D$ either belongs to $D$ or has odd order, then $B$ cannot be a product of $\ell$ isomorphic non-abelian simple groups with $\ell>1$. So $\ell=1$ and $L=T$. However, by \cite[Proposition~1.2]{PraegerSY}, there is no non-abelian simple group with a conjugaucy class $D$ of transpositions such that the product of any two distinct members of $D$ either belongs to $D$ or has order $n$ for a fixed integer $n$. A contradiction occurs. 


\medskip
\f{\bf Step~2.}\ $G$ is not of type {\rm CD}.

Suppose on the contrary that $G$ is quasiprimitive on $V$ of type {\rm CD}. Then we may let $V=\D^k$ with $k\geq 2$, and let $N=T^{k\ell}\unlhd G\leq H\wr S_k$, where $H$ is a quasiprimitive group of type SD on $\D$ with a unique minimal normal subgroup $M=T^\ell$. Furthermore, we may let
\[\D=\{(t_1,\dots,t_{\ell-1}, 1)\mid t_i\in T, i=1,2,\dots,\ell-1\},\]
and the action of $M$ on $\D$ is given by
\[(s_1,\dots,s_{\ell-1},s_\ell):\ (t_1,\dots,t_{\ell-1}, 1)\mapsto (s_\ell^{-1}t_1s_1,\dots,s_\ell^{-1}t_{\ell-1}s_{\ell-1}, 1).\]
Note that $\D$ is a normal subgroup of $M$ which is isomorphic to $T^{\ell-1}$ and acts regularly on $\D$ by right multiplication. Then $\D^{k}$ acts regularly on $V=\D^k$  by right multiplication. So $\G=\Cay(\D^k, S)$ is a Cayley graph of $\D^k$.

Take $\d=(1,\dots,1,1)\in\D$. Then $M_\d=\{(t,\dots,t)\mid t\in T\}$.
Now take $u=(\d,\dots,\d)$. Then $u$ is the identity element of $\D^k=V$ and $N_u=M_\d^k\cong T^k$. Since $\D^k$ acts regularly on itself by right multiplication, one has $N=\D^k: N_u$. Since $\G$ is connected, the neighbourhood $\G(u)=S$ is a generating set of $\D^k$ and $N_u\leq\Aut(\D^k,S)$. This implies that $N_u$ acts faithfully on $\G(u)$. Since $N\unlhd G$, one has $N_u\unlhd G_u$. 
We obtain the following observation.

\medskip\f{\bf Claim~1.}\ $V=\lg\G(u)\rg$, $N_u=M_\d^k\cong T^k$, $N_u\leq\Aut(\D,S)$ and $N_u\unlhd G_u$.
\medskip


Recall that $\G$ is of girth $3$ and is neither complete nor complete multipartite. Then Theorem~\ref{local-pty}~(iv) or (v) happens.

Suppose first that Theorem~\ref{local-pty}~(iv) happens. Then $[\G(u)]$ is disconnected, and $G_u^{\G(u)}$ is imprimitive of rank $3$ and $[\G(u)]\cong r\K_b$ for some integers $r,b\geq 2$. Let $\G(u)=\bigcup_{i=1}^r\Om_i$ such that $[\Om_i]\cong \K_b$. Then $\Sigma=\{\Om_i\ |\ 1\leq i\leq r\}$ is a system of blocks of imprimitivity of $G_{u}$ on $\G(u)$. Since $\G$ is $(G,3)$-CH, $(G_{u})_{\Om_i}$ is $2$-transitive on $\Om_i$ for all $i$, and $G_{u}$ acts $2$-transitively on $\Sigma$. As $N_u=M_\d^k$ with $M_\d\cong T$ and $k\geq 2$, by Proposition~\ref{Burnside}, $N_u$ is not faithful on $\Sigma$. Let $L=\{(g,1,\dots,1)\in T^{k\ell}\mid g\in M_\d\}$. Then $L\cong M_\d$ and $L\leq N_u$.

Without loss of generality, we may assume that $L$ acts trivially on $\Sigma$. Since $L\leq N_u\leq\Aut(\D,S)$, $L$ acts faithfully on $S$. So $L$ acts non-trivially on at least one $\Om_i$. Without loss of generality, assume that $L$ acts non-trivially on $\Om_1$. Since $L\cong T$ is non-abelian simple, $L$ acts faithfully on $\Om_1$. Then
\[L\cong L(G_u)_{(\Om_1)}/(G_u)_{(\Om_1)}\unlhd (N_u)_{\Om_1}(G_u)_{(\Om_1)}/(G_u)_{(\Om_1)}.\]
By Claim~1, we have $N_u\unlhd G_u$. It follows that $(N_u)_{\Om_1}(G_u)_{(\Om_1)}/(G_u)_{(\Om_1)}$ is a non-abelian normal subgroup of $(G_u)_{\Om_1}/(G_u)_{(\Om_1)}$. Since $(G_u)_{\Om_1}$ acts $2$-transitively on $\Om_1$, by Proposition~\ref{Burnside}, $\soc((G_u)_{\Om_1}/(G_u)_{(\Om_1)})$ is non-abelian simple. It implies that $(G_u)_{\Om_1}/(G_u)_{(\Om_1)}$ is almost simple and
\[\soc((G_u)_{\Om_1}/(G_u)_{(\Om_1)})\leq (N_u)_{\Om_1}(G_u)_{(\Om_1)}/(G_u)_{(\Om_1)}\leq (G_u)_{\Om_1}/(G_u)_{(\Om_1)}.\]
As $L(G_u)_{(\Om_1)}/(G_u)_{(\Om_1)}\unlhd (N_u)_{\Om_1}(G_u)_{(\Om_1)}/(G_u)_{(\Om_1)}$, one has
\[L\cong L(G_u)_{(\Om_1)}/(G_u)_{(\Om_1)}=\soc((G_u)_{\Om_1}/(G_u)_{(\Om_1)}).\]

Since $L=\{(g,1,\dots,1)\in T^{k\ell}\mid g\in M_\d\}$ acts non-trivially on $\Om_1$, we may take $v=(t_{1},\dots, t_{\ell-1},1,\ldots)\in \Om_1$ such that $t_i\neq 1$ for some $1\leq i\leq \ell-1$. Then
\[L_v\leq \{((t,\dots,t),1\dots,1)\in T^{k\ell}\mid (t,\dots,t)\in M_\d, tt_i=t_it\}\cong C_T(t_i).\]
It follows that $L_v$ centralises $((t_i,\dots,t_i),1,\dots,1)\in L$, where $(t_i,\dots,t_i)\in M_\d$.
Again, since $(G_u)_{\Om_1}$ is $2$-transitive on $\Om_1$, by Proposition~\ref{cor-2-tran-socle}~(2)-(3), $L\cong L(G_u)_{(\Om_1)}/(G_u)_{(\Om_1)}$ is primitive on $\Om_1$. It follows that $L_v$ is maximal in $L$. Since $L\cong T$ is non-abelian simple, $L_v$ is equal to the centraliser of $((t_i,\dots,t_i),1,\dots,1)$ in $L$. However, this is not possible by Proposition~\ref{cor-2-tran-socle}~(1).

Suppose now Theorem~\ref{local-pty}~(v) happens. Then $[\G(\a)]$ is connected, and $G_\a^{\G(\a)}$ is primitive of rank $3$. By Claim~1, we see that $N_u\cong T^k (k\geq 2)$ is isomorphic to a normal subgroup of $G_u^{\G(u)}$. Then $G_u^{\G(u)}\leq H \wr \mz_2$, where $H$ is an almost simple $2$-transitive permutation group of degree $n_0$, and $|\G(u)|=n_0^2$ (see, for example, \cite[p.165]{Liebeck-Saxl}). So $\soc(G_u^{\G(u)})=\soc(H)^2=N_u$. It follows that $k=2$ and $N_u$ is transitive on $\G(u)$. So $\G(u)=v^{N_u}$ for some $v\in\G(u)$. 
Let $$v=(t_{1},\ldots, t_{\ell-1},1,t_{\ell+1},\ldots,t_{2\ell-1},1)\in\G(u).$$
Then $\G(u)=v^{N_u}=v_1^{M_\d}\times v_2^{M_\d}$, where $v_1=(t_{1},\ldots, t_{\ell-1},1)\in\D$ and $v_2=(t_{\ell+1}, $ $ \ldots, $ $t_{2\ell-1},$ $ 1)\in \D$. So, $|v_1^{M_\d}|=|v_2^{M_\d}|=n_0$, and $H$ is $2$-transitive on $v_i^{M_\d}$ with socle $M_\d$. By Proposition~\ref{cor-2-tran-socle}~(2)-(3), $M_\d$ is primitive on $v_i^{M_\d}$, and so $(M_\d)_{v_i}$ is maximal in $M_\d$. Since $\lg\G(u)\rg=\D^2$ and $M_\d$ is the full diagonal subgroup of $M$, all $t_i$'s are pair-wise distinct and none of them are identity. For each $i\in\{1,\dots,\ell-1\}$, considering the stabiliser $(M_\d)_{v_1}$ of $v_1$ in $M_\d$, we have $(M_\d)_{v_1}\cong C_T(t_{i})$. This is not possible by Proposition~\ref{cor-2-tran-socle}~(1). \hfill\qed


\section{Quasiprimitive $3$-CH graphs of type SD}\label{sec:part2}
In this section, we shall give a classification of quasiprimitive $3$-CH graphs of type SD. In particular, we prove part (2) of Theorem~\ref{th-quasiprimitive}. Before proceeding, we first describe the quasiprimitive permutation groups of type SD (see \cite[Sections~6 \& 12]{Praeger}). Let $T$ be a finite non-abelian simple group and let $\ell\geq 2$ be an integer. A quasiprimitive permutation group is a subgroup of the group $$W:=\{(\a_1,\ldots,\a_\ell)\s\ |\ \a_i\in\Aut(T), \s\in S_\ell, \a_i\Inn(T)=\a_j\Inn(T), 1\leq i,j\leq \ell\},$$
where $\s^{-1}(\a_1,\ldots,\a_\ell)\s=(\a_{1^{\s^{-1}}},\ldots,\a_{\ell^{\s^{-1}}})$. Further, $\soc(W)=\{(t_1,\ldots,t_\ell)\mid t_i\in \Inn(T), 1\leq i\leq \ell\}$. For each $1\leq i\leq \ell$, let $T_i:=\{(1,\ldots,1,\mathop t\limits^{i},1,\ldots,1)\ |\ t\in T\}$ and let $V=T_1\times\cdots\times T_{\ell-1}$  (which we identify with $\Inn(T_1)\times\cdots\times \Inn(T)^{\ell-1}$). Then a primitive action of $W$ on $V$ is defined by
$$\begin{array}{rl}
(\a_1, \ldots, \a_\ell):& (t_1, \ldots, t_{\ell-1},1)\mapsto (\a_\ell^{-1}t_1\a_1, \ldots, \a_\ell^{-1}t_{\ell-1}\a_{\ell-1}, 1),\ {\rm and}\\

\s:& (t_1, \ldots, t_{\ell-1},t_\ell)\mapsto (t_{\ell^{\s^{-1}}}^{-1}t_{1^{\s^{-1}}}, \ldots, t_{\ell^{\s^{-1}}}^{-1}t_{(\ell-1)^{\s^{-1}}},1),\ {\rm where}\ t_\ell=1,
\end{array}
$$
for $(\a_1, \ldots, \a_\ell), \s\in W$. For $u=(1,\dots,1)\in V=T^{\ell-1}$, we have $W_u=A\times S_\ell$, where $A=\{(a,\ldots,a)\ |\ a\in\Aut(T)\}$. The condition for a subgroup $G$ of $W$ which contains $\soc(W)$ is quasiprimitive on $V=T^{\ell-1}$ is that $G$ acts transitively by conjugation on the $\ell$ simple director factors of $\soc(W)$.

\begin{defi}\label{def:sd-graph}
{\rm Under the notation above, let $\G=\Cay(V, S)$ be a connected Cayley graph of $V=T_1\times\cdots\times T_{\ell-1}$, where $T_i:=\{(1,\ldots,1,\mathop t\limits^{i},1,\ldots,1)\ |\ t\in T\}$ with $1\leq i\leq \ell-1$. We say that $\G$ is a {\em $G$-quasiprimitive graph of type SD} with $G\leq\Aut(\G)$ if $G$ acts transitively by conjugation on the $\ell$ simple director factors of $\soc(W)$. }
\end{defi}

\begin{lem}\label{lem:cd-general}
Let $\G=\Cay(V,S)$ be a connected $G$-quasiprimitive graph of type SD. Let $N=\soc(G)=T^\ell$ with $\ell\geq2$ and let $u=(1,\dots,1)$ be the identity of $V$. Let $W=N_{{\rm Sym}(V)}(N)$. For $w\in W_u$, if $w$ fixes $S$ setwise, then $w\in\Aut(\G)$.
\end{lem}

\f\demo By Definition~\ref{def:sd-graph}, $V=T_1\times\cdots\times T_{\ell-1}$, where $T_i:=\{(1,\ldots,1,\mathop t\limits^{i},1,\ldots,1)\ |\ t\in T\}$ with $1\leq i\leq \ell-1$. Further, $N$ is transitive on $V$ with stabiliser $N_u=\{(t,\dots,t)\mid t\in T\}$.

Assume that $w\in W_u$ fixes $S$ setwise. Assume that $\{g, h\}$ is an edge of $\G$. Then $h=sg$ for some $s\in S$. Then $s^w=s'$ for some $s'\in S$ as $S^w=S$. Since $N\unlhd W$ and $w\in W_u$, one has $h^wN_u=h^wN_u^w=(hN_u)^w=(sgN_u)^w=s^wg^wN_u^w=s^wg^wN_u$. It implies that $(h^w)^{-1}s^wg^w\in N_u$. As $s,g,h\in V$, one has $h^w,s^w,g^w\in V^w$ and hence $(h^w)^{-1}s^wg^w\in V^w$. Since $V\unlhd N$, one has $V^w\unlhd N^w=N$, and so $V^w\cap N_u=1$ as $N_u\cong T$ is the full diagonal subgroup of $N$. Consequently, we have $(h^w)^{-1}s^wg^w\in V^w\cap N_u=1$, and hence $h^w=s^wg^w$. This implies that $\{g,h\}^w=\{g,sg\}^w=\{g^w,s^wg^w\}=\{g^w,s'g^w\}$, and $\{g, h\}^w$ is also an edge of $\G$. Conversely, if $\{g,h\}^w$ is an edge of $\G$, then since $S^{w^{-1}}=S$, the argument above implies that $\{g,h\}$ is also an edge of $\G$. Thus, $w\in\Aut(\G)$.\hfill\qed

Next we give four examples of $(G, 3)$-CH and $G$-quasiprimitive graphs of type SD.

\begin{exam}\label{exam:sd-a5-1}
{\rm Let $T=A_5$. Then $T$ has two involutions $a,b$ and an element $c$ of order $3$ such that $\lg a, b\rg\cong\mz_2\times\mz_2$ and $a^c=b, b^c=ab$. Furthermore, $T$ has an automorphism $\a$ such that $a^\a=b,b^\a=a,c^\a=c^{-1}$.

Let $N=T^3$, let $V=\{(t_1,t_2,1)\mid t_1,t_2\in T\}$ and let $S=\{(a, b, 1)^{\bf t}\mid {\bf t}\in\{(t,t,t)\mid t\in T\}\}$. Define $\G_{A_5^3}=\Cay(V, S)$. It is easy to see that
\[S=\bigcup_{g\in Q}\{(a,b,1), (b,ab,1), (ab,a,1)\}^{(g,g,g)}, \]
where $Q$ is a Sylow $5$-subgroup of $T$. So $|S|=15$.

Let $\s=(1, 2, 3), \pi=(1,2)\in S_3$ and $x=(\a,\a,\a)$.  For any $(a^t, b^t, 1)\in S$, we have
\[\begin{array}{l}
(a^t, b^t, 1)^\s=(b^t, a^tb^t, 1)=(a^{ct}, b^{ct}, 1)=(a,b,1)^{(ct,ct,ct)}\in S,\\
(a^t,b^t,1)^{\a\pi}=((a^t)^{\a}, (b^t)^\a, 1)^\pi=((b^t)^{\a}, (a^t)^\a, 1)=(a,b,1)^{(\a^{-1}t\a, \a^{-1}t\a,\a^{-1}t\a)}\in S,\\
(a,b,1)^{\pi}=(b,a,1)\notin S.
\end{array}
\]
By Lemma~\ref{lem:cd-general} we see that $\s, \a\pi\in\Aut(\G_{A_5^3})$. Since $(a,b,1)^{\pi}=(b,a,1)\notin S$, we have $\pi\notin \Aut(\G_{A_5^3})$. Let $G=N_{{\rm Sym}(V)}(N)\cap \Aut(\G_{A_5^3})$. Then $G=N: (\lg\s\rg:\lg \a\pi\rg)$, and $\G_{A_5^3}$ is a $G$-quasiprimitive graph of type SD. Let $u=(1,1,1)$. Then $G_u=(N_u\times\lg\s\rg)\rtimes\lg \a\pi\rg\cong (A_5\times \mz_3): \mz_2$ and $N_u=\{(t,t,t)\mid t\in T\}$ is transitive on $S$ as $S=\{(a, b, 1)^g\mid g\in N_u\}$. It is easy to see that $\soc(G_u)=N_u\times\lg\s\rg$ which is faithful on $S$. Then $G_u$ is faithful on $S$, and so we may view $G_u$ as a transitive permutation group on $S$. Since $N_u$ is transitive on $S$ with point stabiliser a Sylow $2$-subgroup of $N_u$, by Magma~\cite{BCP}, we see that $G_u$ is the normaliser of $N_u$ in ${\rm Sym}(S)$ and $G_u$ has rank $3$. It is easy to see that $\G$ has girth $3$. By Proposition~\ref{prop-local-action}, $\G_{A_5^3}$ is $(G,3)$-CH.}
\end{exam}

\begin{exam}\label{exam:sd-a5-2}
{\rm Let $T=A_5$. Then $T$ has two involutions $a,b$ and an element $c$ of order $3$ such that $\lg a, b\rg\cong\mz_2\times\mz_2$ and $a^c=b, b^c=ab$. Furthermore, $T$ has an automorphism $\a$ such that $a^\a=b,b^\a=a,c^\a=c^{-1}$.

Let $N=T^4$, let $V=\{(t_1,t_2, t_3, 1)\mid t_1,t_2,t_3\in T\}$ and let $S=\{(a, b, ab, 1)^{\bf t}\mid {\bf t}\in\{(t,t,t,t)\mid t\in T\}\}$. Define $\G_{A_5^4}=\Cay(V, S)$. It is easy to see that
\[S=\bigcup_{g\in Q}\{(a,b,ab, 1), (b,ab,a,1), (ab,a,b,1)\}^{(g,g,g,g)}, \]
where $Q$ is a Sylow $5$-subgroup of $T$. So $|S|=15$.

Let $\s=(1, 2, 3), \d_1=(1,4)(2,3), \d_2=(1,3)(2,4), \pi=(1,2)\in S_4$ and $x=(\a,\a,\a,\a)$.  For any $(a^t, b^t, a^tb^t, 1)\in S$, we have
\[\begin{array}{l}
(a^t, b^t, a^tb^t, 1)^\s=(a^tb^t, a^t, b^t, 1)=(b^{ct}, (ab)^{ct}, a^{ct}, 1)=(b, ab, a, 1)^{(ct,ct,ct)}\in S,\\
(a^t, b^t, a^tb^t, 1)^{\d_1}=(a^t, a^t(a^tb^t), a^tb^t, 1)=(a^t, b^t, a^tb^t,1)\in S,\\
(a^t, b^t, a^tb^t, 1)^{\d_2}=(b^t(a^tb^t), b^t, b^ta^t, 1)=(a^t, b^t, a^tb^t,1)\in S,\\
(a^t,b^t, a^tb^t, 1)^{\a\pi}
=((b^t)^{\a}, (a^t)^\a, (a^tb^t)^\a, 1)=(a,b, ab, 1)^{(\a^{-1}t\a, \a^{-1}t\a,\a^{-1}t\a)}\in S,\\
(a,b,ab,1)^{\pi}=(b,a,ab, 1)\notin S.
\end{array}
\]
By Lemma~\ref{lem:cd-general} we see that $\s, \d_1,\d_2, \a\pi\in\Aut(\G_{A_5^4})$, and $\d_1,\d_2$ fixes $S$ pointwise. Since $(a,b,ab,1)^{\pi}=(b,a,ab, 1)\notin S$, we have $\pi\notin \Aut(\G_{A_5^4})$. Let $G=N_{{\rm Sym}(V)}(N)\cap \Aut(\G_{A_5^4})$. Then $G=N: (\lg\s,\d_1,\d_2\rg:\lg \a\pi\rg)$, and $\G_{A_5^4}$ is a $G$-quasiprimitive graph of type SD. Let $u=(1,1,1,1)$. Then $G_u=(N_u\times\lg\s,\d_1,\d_2\rg)\rtimes\lg \a\pi\rg\cong (A_5\times A_4): \mz_2$ and $N_u=\{(t,t,t,t)\mid t\in T\}$ is transitive on $S$ as $S=\{(a, b, ab, 1)^g\mid g\in N_u\}$. It is easy to see that the kernel of $G_u$ on $S$ is $\lg \d_1,\d_2\rg$. So $G_u^S\cong (A_5\times \mz_3): \mz_2$. Since $N_u$ is transitive and faithful on $S$ with point stabiliser a Sylow $2$-subgroup of $N_u$, by Magma~\cite{BCP}, we see that $G_u^S$ is the normaliser of $N_u$ in ${\rm Sym}(S)$ and $G_u^S$ has rank $3$. It is easy to see that $\G_{A_5^4}$ has girth $3$. By Proposition~\ref{prop-local-action}, $\G_{A_5^4}$ is $(G,3)$-CH.}
\end{exam}

\begin{exam}\label{exam:sd-u34-1}
{\rm Let $T=\PSU(3,4)$. Then the center $Z$ of a Sylow $2$-subgroup of $T$ is isomorphic to $\mz_2\times\mz_2$. Set $Z=\lg a,b\rg$. Then $N_T(Z)$ has an element $c$ of order $3$ such that $a^c=b, b^c=ab$. Furthermore, $T$ has an automorphism $\a$ of order $4$ such that $a^\a=b, b^\a=a,c^\a=c^{-1}$.

Let $N=T^3$, let $V=\{(t_1,t_2, 1)\mid t_1,t_2\in T\}$ and let $S=\{(a, b, 1)^{\bf t}\mid {\bf t}\in\{(t,t,t)\mid t\in T\}\}$. Define $\G_{U(3,4)^3}=\Cay(V, S)$. Note that $C_G(a)=C_G(Z)$ has order $2^6\cdot 5$. It follows that the conjugate class $a^T$ of $T$ containing $a$ has length $165$. Clearly, $a, b, ab\in a^T$ and $Z=\lg a,b\rg$. Set $\Om_1=\{(a,b,1), (b,ab,1), (ab,a,1)\}$. Then the setwise stabiliser of $\Om_1$ in the full diagonal subgroup $\{(t,t,t)\mid t\in T\}$ is isomorphic to $[4^3]:\mz_{15}$ which is maximal in $T$.
So $S$ can be partitioned into $65$ subsets of cardinality $3$:
\[S=\bigcup_{t\in T}\{(a,b,1), (b,ab,1), (ab,a,1)\}^{(t,t,t)}.\]
It implies that $|S|=195$.

Let $\s=(1, 2, 3), \pi=(1,2)\in S_3$ and $x=(\a,\a,\a)$.  For any $(a^t, b^t, 1)\in S$, we have
\[\begin{array}{l}
(a^t, b^t, 1)^\s=(b^t, a^tb^t, 1)=(a^{ct}, b^{ct}, 1)=(a,b,1)^{(ct,ct,ct)}\in S,\\
(a^t,b^t,1)^{\a\pi}=((a^t)^{\a}, (b^t)^\a, 1)^\pi=((b^t)^{\a}, (a^t)^\a, 1)=(a,b,1)^{(\a^{-1}t\a, \a^{-1}t\a,\a^{-1}t\a)}\in S,\\
(a,b,1)^{\pi}=(b,a,1)\notin S.
\end{array}
\]
By Lemma~\ref{lem:cd-general} we see that $\s,\a\pi\in\Aut(\G_{U(3,4)^4})$. Since $(a,b,ab,1)^{\pi}=(b,a,ab, 1)\notin S$, we have $\pi\notin \Aut(\G_{U(3,4)^4})$. Let $G=N_{{\rm Sym}(V)}(N)\cap \Aut(\G_{U(3,4)^4})$. Then $G=N: (\lg\s\rg:\lg \a\pi\rg)$, and $\G_{U(3,4)^4}$ is a $G$-quasiprimitive graph of type SD. Let $u=(1,1,1)$. Then $G_u=(N_u\times\lg\s\rg)\rtimes\lg \a\pi\rg\cong (\PSU(3,4)\times \mz_3): \mz_4$ and $N_u=\{(t,t,t)\mid t\in T\}$ is transitive on $S$ as $S=\{(a, b, 1)^g\mid g\in N_u\}$. It is easy to see that $G_u$ is faithful on $S$. So $G_u^S\cong G_u\cong (\PSU(3,4)\times \mz_3): \mz_4$. Since $N_u$ is transitive and faithful on $S$ with point stabiliser a Sylow $2$-subgroup of $N_u$, by Magma~\cite{BCP}, we see that $G_u^S$ is the normaliser of $N_u$ in ${\rm Sym}(S)$ and $G_u^S$ has rank $3$. It is easy to see that $\G_{U(3,4)^3}$ has girth $3$. By Proposition~\ref{prop-local-action}, $\G_{U(3,4)^3}$ is $(G,3)$-CH.}
\end{exam}

\begin{exam}\label{exam:sd-u34-2}
{\rm Let $T=\PSU(3,4)$. Then the center $Z$ of a Sylow $2$-subgroup of $T$ is isomorphic to $\mz_2\times\mz_2$. Set $Z=\lg a,b\rg$. Then $N_T(Z)$ has an element $c$ of order $3$ such that $a^c=b, b^c=ab$. Furthermore, $T$ has an automorphism $\a$ of order $4$ such that $a^\a=b, b^\a=a,c^\a=c^{-1}$.

Let $N=T^4$, let $V=\{(t_1,t_2, t_3, 1)\mid t_1,t_2,t_3\in T\}$ and let $S=\{(a, b, ab, 1)^{\bf t}\mid {\bf t}\in\{(t,t,t,t)\mid t\in T\}\}$. Define $\G_{U(3,4)^4}=\Cay(V, S)$. Note that $C_G(a)=C_G(Z)$ has order $2^6\cdot 5$. It follows that the conjugate class $a^T$ of $T$ containing $a$ has length $165$. Clearly, $a, b, ab\in a^T$ and $Z=\lg a,b\rg$. Set $\Om_1=\{(a,b,ab,1), (b,ab,a,1), (ab,a,b,1)\}$. Then the setwise stabiliser of $\Om_1$ in the full diagonal subgroup $\{(t,t,t,t)\mid t\in T\}$ is isomorphic to $[4^3]:\mz_{15}$ which is maximal in $T$.
So $S$ can be partitioned into $65$ subsets of cardinality $3$:
\[S=\bigcup_{t\in T}\{(a,b,ab,1), (b,ab,a,1), (ab,a,b,1)\}^{(t,t,t,t)}.\]
It implies that $|S|=195$.

Let $\s=(1, 2, 3), \d_1=(1,4)(2,3), \d_2=(1,3)(2,4), \pi=(1,2)\in S_4$ and $x=(\a,\a,\a,\a)$.  For any $(a^t, b^t, a^tb^t, 1)\in S$, we have
\[\begin{array}{l}
(a^t, b^t, a^tb^t, 1)^\s=(a^tb^t, a^t, b^t, 1)=(b^{ct}, (ab)^{ct}, a^{ct}, 1)=(b, ab, a, 1)^{(ct,ct,ct)}\in S,\\
(a^t, b^t, a^tb^t, 1)^{\d_1}=(a^t, a^t(a^tb^t), a^tb^t, 1)=(a^t, b^t, a^tb^t,1)\in S,\\
(a^t, b^t, a^tb^t, 1)^{\d_2}=(b^t(a^tb^t), b^t, b^ta^t, 1)=(a^t, b^t, a^tb^t,1)\in S,\\
(a^t,b^t, a^tb^t, 1)^{\a\pi}
=((b^t)^{\a}, (a^t)^\a, (a^tb^t)^\a, 1)=(a,b, ab, 1)^{(\a^{-1}t\a, \a^{-1}t\a,\a^{-1}t\a)}\in S,\\
(a,b,ab,1)^{\pi}=(b,a,ab, 1)\notin S.
\end{array}
\]
By Lemma~\ref{lem:cd-general} we see that $\s, \d_1,\d_2, \a\pi\in\Aut(\G_{U(3,4)^4})$, and $\d_1,\d_2$ fixes $S$ pointwise. Since $(a,b,ab,1)^{\pi}=(b,a,ab, 1)\notin S$, we have $\pi\notin \Aut(\G_{U(3,4)^4})$. Let $G=N_{{\rm Sym}(V)}(N)\cap \Aut(\G_{U(3,4)^4})$. Then $G=N: (\lg\s,\d_1,\d_2\rg:\lg \a\pi\rg)$, and $\G_{U(3,4)^4}$ is a $G$-quasiprimitive graph of type SD. Let $u=(1,1,1,1)$. Then $G_u=(N_u\times\lg\s,\d_1,\d_2\rg)\rtimes\lg \a\pi\rg\cong (\PSU(3,4)\times A_4): \mz_4$ and $N_u=\{(t,t,t,t)\mid t\in T\}$ is transitive on $S$ as $S=\{(a, b, ab, 1)^g\mid g\in N_u\}$. It is easy to see that the kernel of $G_u$ on $S$ is $\lg \d_1,\d_2\rg$. So $G_u^S\cong (\PSU(3,4)\times \mz_3): \mz_4$. Since $N_u$ is transitive and faithful on $S$ with point stabiliser a Sylow $2$-subgroup of $N_u$, by Magma~\cite{BCP}, we see that $G_u^S$ is the normaliser of $N_u$ in ${\rm Sym}(S)$ and $G_u^S$ has rank $3$. It is easy to see that $\G_{U(3,4)^4}$ has girth $3$. By Proposition~\ref{prop-local-action}, $\G_{U(3,4)^4}$ is $(G,3)$-CH.}
\end{exam}

The following is the main theorem of this section which is just part (2) of Theorem~\ref{th-quasiprimitive}.

\begin{theorem}\label{SD}
Let $\Ga=(V,E)$ be a connected $(G,3)$-CH graph of girth $3$, which is neither complete nor complete multipartite. If $G$ is quasiprimitive on $V$ of type SD, then $\G\cong\G_{A_5^3}, \G_{A_5^4}, \G_{U(3,4)^3}$ or $\G_{U(3,4)^4}$ $($as given in {\rm Examples~\ref{exam:sd-a5-1}--\ref{exam:sd-u34-2})}.
\end{theorem}

\f\demo Suppose that $G$ is quasiprimitive on $V$ of type {\rm SD}.
Let $N=\soc(G)=T^\ell$ with $\ell\geq 2$ and $T$ a nonabelian simple group. Then $N$ is the unique minimal normal subgroup of $G$. Furthermore, we may let
\[V=\{(t_1,\dots,t_{\ell-1}, 1)\mid t_i\in T, i=1,2,\dots,\ell-1\}\},\]
and the action of $N$ on $V$ is given by
\[(s_1,\dots,s_{\ell-1},s_\ell):\ (t_1,\dots,t_{\ell-1}, 1)\mapsto (s_\ell^{-1}t_1s_1,\dots,s_\ell^{-1}t_{\ell-1}s_{\ell-1}, 1).\]
Note that $V$ is a normal subgroup of $N$ which is isomorphic to $T^{\ell-1}$ and acts regularly on $V$ by right multiplication. So $\G=\Cay(V, S)$ is a Cayley graph of $V$. Take $u=(1,\dots,1,1)\in V$. Then $u$ is the identity element of $V$ and $N_u=\{(t,\dots,t)\mid t\in T\}$ is the full diagonal subgroup of $N$. So $N_u\cong T$ and $N=V: N_u$. Since $\G$ is connected, the neighbourhood $S=\G(u)$ generates $V$. It follows that $N_u$ acts faithfully on $\G(u)$ as $N_u\leq\Aut(V,\G(u))$. We have the following observation.\medskip

\f{\bf Claim~1.}\ $V=\lg\G(u)\rg$, $N_u\cong T$, $N_u\leq\Aut(V,\G(u))$ and $N_u\unlhd G_u$.
\medskip




Next we will prove the following claim.

\medskip
\f{\bf Claim~2.}\ If $N_u^{\G(u)}$ is transitive, then either $\ell>2$ or $|w|>2$ for any $w\in \G(u)$.\medskip

To prove this claim, we suppose that $N_u^{\G(u)}$ is transitive, and $\ell=2$ and $|w|=2$. Then $V\cong T$ and $V\leq V:\Aut(V,\G(u))\leq V: \Aut(V)$. So $V: \Aut(V, \G(u))$ is quasiprimitive on $V$ of type HS. By Theorem~\ref{th:hs-hc-cd}, $\G$ is not $(V: \Aut(V, \G(u)), 3)$-CH. This implies that $G_u\nleq\Aut(V)$.

Let $B=\{(a, a)\ |\ a\in\Aut(T)\}$ and take $1\neq\s\in S_2$. Then $B$ acts on $V$ by
\[(a, a):\ (t_1, 1)\mapsto (t_1^a,1),\]
and $\s$ acts on $V$ by
\[\s: (t_1, 1)\mapsto (t_1^{-1},1).\]
So $N\leq G\leq W=N_{{\rm Sym}(V)}(N)=VW_u$, where $W_u=B\times\lg \s\rg$. Clearly, $B$ normalises $V$, so $B\leq\Aut(V)$. Since  $G_u\nleq\Aut(V)$, one has $G_u\nleq B$. For any $g\in G_u\setminus B$, we have $g=(a, a)\s$ for some $a\in\Aut(T)$. Take $w=(t_1,1)\in S$. Since $w$ has order $2$, $t_1$ is an involution, and so $w^\s=(t_1, 1)^\s=(t_1^{-1}, 1)=(t_1, 1)=w$. This implies that $\s$ fixes $S$ pointwise. Since $g\in G_u$, one has $w^g\in \G(u)$ and so $w^{(a, a)}\in \G(u)$. By the arbitrariness of $w$, we have $\G(u)^{(a, a)}=\G(u)$, and hence $(a,a)\in \Aut(V, \G(u))$. It follows that $G_u\leq \Aut(V,\G(u))\times\lg \s\rg$. Since $G_u^{\G(u)}$ is a transitive permutation group of rank $3$,
$(\Aut(V,\G(u))\times\lg \s\rg)^{\G(u)}$ is also a transitive permutation group of rank $3$. As $\s$ fixes $\G(u)$ pointwise, $\Aut(V,\G(u))^{\G(u)}$ is a transitive permutation group of rank 3, and since $\G$ is assumed to have girth $3$, by Proposition~\ref{prop-local-action}, $\G$ is $(V: \Aut(H,\G(u)), 3)$-CH, a contradiction.\medskip

Recall that $\G$ is of girth $3$ and is neither complete nor complete multipartite. Then Theorem~\ref{local-pty}~(iv) or (v) happens.

We first suppose that Theorem~\ref{local-pty}~(v) happens. Then $[\G(u)]$ is connected and $G_u^{\G(u)}$ is primitive of rank $3$. Then either $\soc(G_u^{\G(u)})$ is abelian, or $G_u^{\G(u)}$ has a unique minimal normal subgroup $R^\ell$, where $\ell=1$ or $2$ and $R$ is a non-abelian simple group (see, for example, \cite[p.165]{Liebeck-Saxl}). By Claim~1, $N_u\cong T$ is isomorphic to a normal subgroup of $G_u^{\G(u)}$. It follows that $G_u^{\G(u)}$ is almost simple with socle $N_u\cong T$. A list of smallest possible primitive rank 3 almost simple groups was given in \cite{Bue-Van} (see also \cite{Devillers}), and by inspecting this list, we conclude that the socle of every primitive rank 3 almost simple group is also primitive. It follows that $N_u$ is primitive on $\G(u)$. So $\G(u)$ consists of involutions and $\G(u)=v^{N_{u}}$ for some $v\in\G(u)$. Let $v=(t_{1},\ldots,t_{\ell-1},1)$. Then each $t_i$ has order $2$. By Claim~1, we have $H=\lg \G(u)\rg=\lg v^{N_u}\rg$. It follows that the $t_{i}$'s are pair-wise distinct and none of them is identity. Since $N_u$ is primitive on $\G(u)$, $N_{uv}$ is maximal in $N_u$. Set $J=\bigcap_{i=1}^{\ell-1}C_T(t_{i})$. Since $N_u\cong T$ is non-abelian simple, one has $N_{uv}\cong J$, and so $J$ is maximal in $T$. Since $v$ is an involution, by Claim~2, we have $\ell>2$, and so $\Z(J)$ contains a subgroup isomorphic to $\mz_2\times\mz_2$. However, by \cite[Theorem~1.1]{Liebeck-Saxl1991}, the center of a point-stabiliser of a primitive permutation group is cyclic. A contradiction occurs.
\medskip



Next we suppose  that Theorem~\ref{local-pty}~(iv) happens. Then $G_u^{\G(u)}$ is imprimitive of rank $3$ and $[\G(u)]\cong r\K_b$ for some integers $r,b\geq 2$. Let $\G(u)=\bigcup_{i=1}^r\Om_i$, where $[\Om_i]\cong \K_b$. Then $\Sigma=\{\Om_i\ |\ 1\leq i\leq r\}$ is a system of blocks of imprimitivity of $G_{u}$ on $\G(u)$. Since $\G$ is $(G,3)$-CH, it follows that $(G_{u})_{\Om_i}$ is $2$-transitive on $\Om_i$ for all $i$, and $G_{u}$ acts $2$-transitively on $\Sigma$. We first have the following claim.\medskip

\f{\bf Claim~3}\ $N_u$ acts primitively and faithfully on $\Sigma$ and $(N_u)_{\Om_1}$ is transitive on $\Om_1$. In particular, $N_u$ is transitive on $\G(u)$.\medskip

Suppose that $N_u$ acts non-faithfully on $\Sigma$. Then $N_u$ would fix each $\Om_i$ setwise because $N_u\cong T$ is non-abelian simple (by Claim~1). As $N_u\leq\Aut(V,\G(u))$ acts faithfully on $\G(u)$, $N_u$ acts faithfully on some $\Om_i$. Note that $(G_{u})_{\Om_i}$ is $2$-transitive on $\Om_i$. As $N_u\unlhd G_u$, $N_u$ is the socle of $(G_{u})_{\Om_i}^{\Om_i}$. By Proposition~\ref{cor-2-tran-socle}~(2)--(3), $N_u$ is primitive on $\Om_i$. Let $v=(t_1, \ldots, t_{\ell-1},1)\in\Om_i$. Then $N_{uv}$ is maximal in $N_u$. Again, considering the stabiliser of $v$ in $N_u$, we have ${\bf t}=(t,\cdots,t)\in N_{uv}$ if and only if $t_i^t=t_i$, that is, $t\in C_T(t_i)$ for all $i$. Since $u\neq v$, at least one of $t_i$'s is not identity. Without loss of generality, assume $t_1\neq 1$. It follows that $N_{uv}\leq C_{N_u}({\bf t}_1)$, where ${\bf t}_1=(t_1,\cdots,t_1)$. Since $N_{uv}$ is maximal in $N_u$ and $N_u$ is non-abelian simple, one has $N_{uv}=C_{N_u}({\bf t}_1)$, which is impossible by Proposition~\ref{cor-2-tran-socle}~(1).

Thus, $N_u$ acts faithfully on $\Sigma$. Since $G_{u}^\Sigma$ is $2$-transitive and $N_u\unlhd G_u$, $N_u$ is the socle of $G_u^\Sigma$. By Proposition~\ref{cor-2-tran-socle}~(2)--(3), $N_u$ is primitive on $\Sigma$ and so $(N_u)_{\Om_1}$ is maximal in $N_u$. Suppose that $(N_u)_{\Om_1}$ fixes ${\Omega_1}$ pointwise. Take $v=(t_1, \ldots, t_{\ell-1},1)\in\Om_1$. Then for any ${\bf t}=(t,\dots, t)\in (N_u)_{\Om_1}$, ${\bf t}$ fixes $v$ and so $t_i^t=t_i$ for all $i$. Again, since $u\neq v$, at least one of $t_i$'s is not identity. Without loss of generality, assume $t_1\neq 1$. It follows that $(N_u)_{\Om_1}\leq C_{N_u}({\bf t}_1)$, where ${\bf t}_1=(t_1, \cdots, t_1)$. Since $(N_u)_{\Om_1}$ is maximal in $N_u$ and $N_u$ is non-abelian simple, one has $(N_u)_{\Om_1}=C_{(N_u)_{\Om_1}}({\bf t}_1)$, which is impossible by Proposition~\ref{cor-2-tran-socle}~(1).

Thus, $(N_u)_{\Om_1}$ does not fix ${\Omega_1}$ pointwise. Since $(G_{u})_{\Om_1}$ is $2$-transitive on $\Om_1$, the normality of $(N_u)_{\Om_1}$ in $(G_{u})_{\Om_1}$ implies that $(N_u)_{\Om_1}$ is transitive on $\Om_1$. It follows that $N_u$ is also transitive on $\G(u)$. \medskip

\f{\bf Claim~4}\ Let $v=(t_1,\ldots, t_{\ell-1},1)\in \Om_1$. Then the $t_i$'s are pair-wise distinct, and none of them is identity.

By Claim~3, we have $\G(u)=v^{N_u}=\{(t_1^t,\ldots,t_{\ell-1}^t,1_T)\ |\ t\in T\}$, and by Claim~1, $\G(u)$ is a generating set of $V$.  As \[V=\{(t_1,\dots,t_{\ell-1}, 1)\mid t_i\in T, i=1,2,\dots,\ell-1\}\}\cong T^{\ell-1},\]
it follows every $t_i$ is non-identity and $t_i\neq t_j$ if $i\neq j$. \medskip

\f{\bf Claim~5}\  Take $v=(t_1, \ldots, t_{\ell-1}, 1)\in \Om_1$. Then $N_u, (N_u)_{\Om_1}, N_{uv}$ and $t_i$ satisfy Table~\ref{table-2} for $i=1,2, \dots, \ell-1$. In particular, for $1\leq i\leq \ell-1$, we have ${\bf t}_i\in Z(N_{uv})$, where ${\bf t}_i=(t_i, \ldots, t_i)$. Furthermore, lines 3-4, 9, 11-12 of Table~\ref{table-2} cannot happen, and $G_u^{\G(u)}$ is not quasiprimitive.

By Claim~4, $t_i\neq 1$ for $1\leq i\leq \ell-1$. Note that $N_{uv}\leq C_{N_u}({\bf t}_i)$, where ${\bf t_i}=(t_i, \ldots, t_i)$ for $1\leq i\leq\ell-1$. So we may apply Theorem~\ref{lem-centralizer} to the $2$-transitive permutation group $G_u^{\Sigma}$, the $2$-transitive action of the point stabiliser $(G_{u})_{\Om_1}$ on $\Om_1$ and the transitive action of $(N_{u})_{\Om_1}$ on $\Om_1$. By Claim~3, $N_u$ is just the socle of $G_u^{\Sigma}$. By Theorem~\ref{lem-centralizer}, $N_u, (N_u)_{\Om_1}, N_{uv}$ and $t_i$ satisfy Table~\ref{table-2}. It follows that for $1\leq i\leq \ell-1$, we have ${\bf t}_i\in Z(N_{uv})$.

If line~3 of Table~\ref{table-2} happens, then $N_u\cong T=A_5$, $|\Om_1|=|(N_u)_{\Om_1}: N_{uv}|=2$, $N_{uv}\cong\mz_5$, $(N_u)_{\Om_1}=D_{10}$ and $r=|N_u: (N_u)_{\Om_1}|=6$. By the above argument, we have $v=(t_1,t_2,\dots,t_{\ell-1},1)$ such that each ${\bf t}_i=(t_i, \ldots, t_i)\in Z(N_{uv})$. By Claim~4, all $t_i$'s are pair-wise distinct and non-identity. This implies that $|v|=5$. Let ${\bf t}=(t,t,\dots,t)\in (N_u)_{\Om_1}$ be an involution. Then ${\bf t}_i^{\bf t}={\bf t}_i^{-1}$ for each $1\leq i\leq\ell-1$. It follows that $v^{-1}=v^{\bf t}\in\Om_1$. So $\Om_1=\{v, v^{-1}\}$. Since $[\Om_1]\cong\K_2$, it follows that $v^2, v^{-2}\in\G(u)$, and so $\{v, v^2\}$ is also an edge of $\G$. This is impossible because $[\G(u)]\cong 6\K_2$.

If line~4 of Table~\ref{table-2} happens, then $N_u\cong T=\PSL(2,8)$, $|\Om_1|=|(N_u)_{\Om_1}: N_{uv}|=2$, $N_{uv}\cong\mz_9$, $(N_u)_{\Om_1}=D_{18}$ and $r=|N_u: (N_u)_{\Om_1}|=28$. It follows that $|\G(u)|=28\cdot2=56$. View $N_u$ as a transitive permutation group on $\G(u)$ with stabiliser $N_{uv}\cong\mz_9$. Then $G_u^{\G(u)}\leq N_{{\rm Sym}(\G(u))}(N_u)$ and $G_u^{\G(u)}$ is of rank $3$. However, by Magma~\cite{BCP}, we see that $N_{{\rm Sym}(\G(u))}(N_u)$ is of rank $4$, a contradiction.

If line~9 of Table~\ref{table-2} happens, then $N_u\cong T=\PSU(3,3)$, $|\Om_1|=|(N_u)_{\Om_1}: N_{uv}|=2$, $N_{uv}\cong[3^3]:\mz_4$, $(N_u)_{\Om_1}=[3^3]:\mz_8$ and $r=|N_u: (N_u)_{\Om_1}|=28$. It follows that $|\G(u)|=28\cdot2=56$. View $N_u$ as a permutation group on $\G(u)$ with point stabiliser $N_{uv}\cong [3^3]:\mz_4$. Then $G_u^{\G(u)}\leq N_{{\rm Sym}(\G(u))}(N_u)$ and $G_u^{\G(u)}$ is of rank $3$. However, by Magma~\cite{BCP}, we see that $N_{{\rm Sym}(\G(u))}(N_u)$ is of rank $4$, a contradiction.

If line~11 of Table~\ref{table-2} happens, then $N_u\cong T=\PSU(3,8)$, $|\Om_1|=|(N_u)_{\Om_1}: N_{uv}|=7$, $N_{uv}\cong[8^3]:\mz_3$, $(N_u)_{\Om_1}=[8^3]:\mz_{21}$ and  $r=|N_u: (N_u)_{\Om_1}|=8^3+1=513$. It follows that $|\G(u)|=513\cdot7=3591$. View $N_u$ as a permutation group on $\G(u)$ with point stabiliser $N_{uv}\cong [8^3]:\mz_3$. Then $G_u^{\G(u)}\leq N_{{\rm Sym}(\G(u))}(N_u)$ and $G_u^{\G(u)}$ is of rank $3$. However, by Magma~\cite{BCP}, we see that $N_{{\rm Sym}(\G(u))}(N_u)$ is of rank $4$, a contradiction.

If line~12 of Table~\ref{table-2} happens, then $N_u\cong T=\PSU(3,32)$, $|\Om_1|=|(N_u)_{\Om_1}: N_{uv}|=31$, $N_{uv}\cong[32^3]:\mz_{11}$, $(N_u)_{\Om_1}=[32^3]:\mz_{341}$  and  $r=|N_u: (N_u)_{\Om_1}|=32^3+1=32769$. It follows that $|\G(u)|=32769\cdot31=1015839$. View $N_u$ as a permutation group on $\G(u)$  with point stabiliser $N_{uv}\cong [32^3]:\mz_3$. Then $G_u^{\G(u)}\leq N_{{\rm Sym}(\G(u))}(N_u)$ and $G_u^{\G(u)}$ is of rank $3$. However, by Magma~\cite{BCP}, we see that $N_{{\rm Sym}(\G(u))}(N_u)$ is of rank $8$, a contradiction.

Finally, suppose that $G_u^{\G(u)}$ is quasiprimitive. Notice that $G_u$ is imprimitive on $\G(u)$. By \cite[Corollary~1.4]{DGLPP}, $G_u^{\G(u)}$ is almost simple, and all quasiprimitive imprimitive permutation groups of rank $3$ are listed in \cite[Table~1]{DGLPP}. For convenience, we give this list as below. Now let $\mathcal{G}$ be a quasiprimitive imprimitive rank 3 permutation group  on a set $\Om$. Then $\mathcal{G}$ preserves a unique non-trivial block system $\mathcal{B}$ on $\Om$ (see \cite[Lemma~3.3]{DGLPP}). Take $B\in\mathcal{B}$ and assume that $|\mathcal{B}|=n$ and $|B|=m$. Then $\mathcal{G}, n, m, \mathcal{G}_B^B$ are as one of the lines of Table~\ref{table-3}.

{\begin{table}[ht]
\caption{Qusiprimitive imprimitive rank 3 groups} \label{table-3}
\center
\newcommand{\tabincell}[2]{\begin{tabular}{@{}#1@{}}#2\end{tabular}}
\begin{tabular}{lllll}
\hline $\mathcal{G}$   &$n$ & $m$ & $\mathcal{G}_B^B$ & extra conditions\\
\hline  $M_{11}$   & $11$ & $2$ & $C_2$& \\
$\mathcal{G}\geq\PSL(2,q)$   & $q+1$ & $2$ & $C_2$& see~\cite[Proposition~5.13]{DGLPP}\\
\hline
$\mathcal{G}\geq\PSL(d,q)$   & $\frac{q^d-1}{q-1}$ & $m$ & ${\rm AGL}(1,m)$ & \tabincell{l}{$d\geq 3$,  $\mathcal{G}$ satisfies the conditions on \\ line 12 of \cite[Table~3]{DGLPP}. In particular,\\
 $m$ prime and $m\mid (q-1)$. } \\
 \hline
 $\PGL(3,4)$   & $21$ & $6$ & $\PSL(2,5)$& \\
 ${\rm P\Gamma L}(3,4)$   & $21$ & $6$ & $\PGL(2,5)$& \\
  $\PSL(3,5)$   & $31$ & $5$ & $S_5$& \\
   $\PSL(5,2)$   & $31$ & $8$ & $A_8$& \\
  ${\rm P\Gamma L}(3,8)$   & $73$ & $28$ & ${\rm Ree}(3)$& \\
   $\PSL(3,2)$   & $7$ & $2$ & $C_2$& \\
   $\PSL(3,3)$   & $13$ & $3$ & $S_3$& \\
\hline
\end{tabular}
\end{table}}

So $G_u^{\G(u)}$ is isomorphic to one of the groups $\mathcal{G}$ in Table~\ref{table-3}. Since $G_u^{\G(u)}$ is almost simple, we have $N_u\cong N_u^{\G(u)}$ which is the socle of $G_u^{\G(u)}$. By the argument above, $N_u, (N_u)_{\Om_1}, N_{uv}$ and $t_i$ satisfy one of lines 1-2, 5-8, 10 of Table~\ref{table-2} with $1\leq i\leq \ell$. By inspecting Table~\ref{table-3}, we see that lines 1-2, 8, 10 of Table~\ref{table-2} do not happen.

If line 5 of Table~\ref{table-2}, then $N_u\cong\PSL(2,q)$, $(N_u)_{\Om_1}=q:\mz_k$ and $N_{uv}=\mz_k$ with $k=\frac{q-1}{(2,q-1)}\neq 2$. So $G_u^{\G(u)}$ has $q+1$ blocks of size $q$ on $\G(u)$. In this case, $G_u^{\G(u)}$ must be isomorphic to the group on line~2 of Table~\ref{table-3}, and then $q=2$, contradicting that $N_u$ is simple.

If line 6 of Table~\ref{table-2}, then $N_u\cong\PSL(d,q)$, $(N_u)_{\Om_1}=q^{d-1}:(\GL(d-1,q)/\mz_k)$ and $N_{uv}=\GL(d-1,q)/\mz_k$ with $d\geq 3, k=(d,q-1)$. So $G_u^{\G(u)}$ has $\frac{q^d-1}{q-1}$ blocks of size $q^{d-1}$ on $\G(u)$. In this case, $G_u^{\G(u)}$ must be isomorphic to the groups on lines~3--10 of Table~\ref{table-3}, but for each of these cases, the block size $m$ is not equal to $q^{d-1}$, a contradiction.

If line 7 of Table~\ref{table-2}, then $N_u\cong\PSL(d,q)$, $(N_u)_{\Om_1}=q^{d-1}: \GL(d-1,q)$ and $N_{uv}=[q^{2d-3}]: \GL(d-2,q)$ with $d\geq 3, (q-1)\mid d$. So $(N_u)_{\Om_1}^{\Om_1}\cong\GL(d-1, q)$ and $G_u^{\G(u)}$ has $\frac{q^d-1}{q-1}$ blocks of size $\frac{q^{d-1}-1}{q-1}$ on $\G(u)$. If $d=3$ and $q=2$, then $|\G(u)|=21$ and $N_u\cong \PSL(3, 2)$. View $N_u$ as a transitive permutation group on $\G(u)$ with point stabiliser $N_{uv}$ a Sylow $2$-subgroup of $N_u$. Then $G_u^{\G(u)}\leq N_{{\rm Sym}(\G(u))}(N_u)$ and $G_u^{\G(u)}$ is of rank $3$. However, by Magma~\cite{BCP}, we see that $N_{{\rm Sym}(\G(u))}(N_u)$ is of rank $4$, a contradiction. Thus, either $d>3$ or $q>2$. It follows that $G_u^{\G(u)}$ is not isomorphic to the group on line $3$ since $(N_u)_{\Om_1}^{\Om_1}\cong\GL(d-1, q)$ which is not isomorphic to ${\rm AGL}(1,m)$. Now as $(q-1)\mid d$, $G_u^{\G(u)}$ must be isomorphic to one of the groups on lines~4--5, 7, 9 of Table~\ref{table-3}, but for each of these cases, the block size $m$ is not equal to $\frac{q^{d-1}-1}{q-1}$, a contradiction.
\medskip

Now we are ready to finish the proof. By Claim~6, $G_u^{\G(u)}$ is not quasiprimitive. Let $M$ be a normal subgroup of $G_u$ which is intransitive and non-trivial on $\G(u)$. Since $[\G(u)]\cong r\K_{b}$, the complement $[\G(u)]^c$ of $[\G(u)]$ is a complete multipartite graph with $r$ parts of size $b$. Since $G_u^{\G(u)}$ is transitive of rank $3$, $G_u$ is arc-transitive on $[\G(u)]^c$. Since $M\unlhd G_u$, each orbit of $M$ on $\G(u)$ does not contains edges of $[\G(u)]^c$, and so each orbit of $M$ on $\G(u)$ induces a complete subgraph of $[\G(u)]$. Recall that $\G(u)=\bigcup_{i=1}^r\Om_i$, where $[\Om_i]\cong \K_b$, and $\Sigma=\{\Om_i\ |\ 1\leq i\leq r\}$ is a system of blocks of imprimitivity of $G_{u}$ on $\G(u)$. Let $B_1$ be an orbit of $M$. Then $|B_1|\neq 1$ and $B_1$ is a block of imprimitivity of $G_u$ on $\G(u)$. Since $[B_1]$ is a complete subgraph of $[\G(u)]$, $B_1$ is contained in some $\Om_i$. Without loss of generality, assume that $B_1\subseteq\Om_1$. As $B_1$ is a block of imprimitivity of $G_u$ on $\G(u)$, $B_1$ is a block of imprimitivity of $(G_u)_{\Om_1}$ on $\Om_1$. Since $(G_u)_{\Om_1}$ is $2$-transitive on $\Om_1$, one has $B_1=\Om_1$. It follows that $\Sigma$ is the set of orbits of $M$ on $\G(u)$.

Since $N_u$ is non-abelian simple and $N_u$ is transitive on $\G(u)$, one has $N_u\cap M=1$, and since $N_u\unlhd G_u$, one has $MN_u=M\times N_u$. So $M\leq C_{G_u}(N_u)$. Since $N_u$ is transitive on $\G(u)$, $M$ is semiregular on $\G(u)$ and so $M$ is regular on each $\Om_i$. Clearly, $M$ is a normal subgroup of $(G_u)_{\Om_1}$. As $(G_u)_{\Om_1}$ is $2$-transitive on $\Om_1$, by Proposition~\ref{Burnside}, $(G_u)_{\Om_1}^{\Om_1}$ has elementary abelian socle. Since $M\leq C_{G_u}(N_u)$ and $M$ is transitive on $\Om_1$, $(N_u)_{\Om_1}$ is semiregular on $\Om_1$. Since $(N_u)_{\Om_1}$ is transitive on $\Om_1$, $(N_u)_{\Om_1}$ is regular on $\Om_1$. Take $v=(t_1,\ldots, t_{\ell-1},1)\in \Om_1$. Then $N_{uv}$ is just the kernel of $(N_u)_{\Om_1}$ acting on $\Om_1$. By Claim~6, $N_u, (N_u)_{\Om_1}, N_{uv}$ and $t_i (1\leq i\leq \ell-1)$ satisfy Table~\ref{table-2}. Since $N_{uv}\unlhd (N_u)_{\Om_1}$, lines 2-4 and 9-12 of Table~\ref{table-2} happen. By Claim~6, lines 2 or $10$ of Table~\ref{table-2} happen.

If line 2 of Table~\ref{table-2} happens, then $N_u\cong T=A_5$, $|\Om_1|=|(N_u)_{\Om_1}: N_{uv}|=3$, $N_{uv}\cong\mz_2\times\mz_2$ and $(N_u)_{\Om_1}=A_4$. By Claim~5, all $t_i$'s are pair-wise distinct, and by Claim~6, for $1\leq i\leq \ell-1$, we have ${\bf t}_i\in Z(N_{uv})$, where ${\bf t}_i=(t_i, \ldots, t_i)$. It follows that $\ell-1\leq |Z(N_{uv})|-1$, and hence $\ell\leq 4$. Since each ${\bf t}_i\in Z(N_{uv})\cong\mz_2\times\mz_2$, one has $|v|=2$. By Claim~3, we see that $\ell=3$ or $4$. If $\ell=3$, then $v=(t_1, t_2, 1)$ where $\lg t_1,t_2\rg\cong\mz_2\times\mz_2$, and then $S=\{(t_1,t_2,1)^{\bf t}\mid {\bf t}=(t,t,t)\in N_u\}$ and $\G\cong\G_{A_5^3}$ (see Example~\ref{exam:sd-a5-1}). If $\ell=4$, then $v=(t_1, t_2, t_3, 1)$ where $\lg t_1,t_2,t_3\rg\cong\mz_2\times\mz_2$, and then $S=\{(t_1,t_2,t_3, 1)^{\bf t}\mid {\bf t}=(t,t,t,t)\in N_u\}$ and $\G\cong\G_{A_5^4}$ (see Example~\ref{exam:sd-a5-2}).

If line~10 of Table~\ref{table-2} happens, then $N_u\cong T=\PSU(3,4)$, $|\Om_1|=|(N_u)_{\Om_1}: N_{uv}|=3$, $N_{uv}\cong[4^3]:\mz_5$ and $(N_u)_{\Om_1}=[4^3]:\mz_{15}$. It follows that $|\G(u)|=195$. By Claim~5, all $t_i$'s are pair-wise distinct, and by Claim~6, for $1\leq i\leq \ell-1$, we have ${\bf t}_i\in Z(N_{uv})\cong\mz_2\times\mz_2$, where ${\bf t}_i=(t_i, \ldots, t_i)$. It follows that $\ell-1\leq |Z(N_{uv})|-1$, and hence $\ell\leq 4$. Since each ${\bf t}_i\in Z(N_{uv})\cong\mz_2\times\mz_2$, one has $|v|=2$. By Claim~3, we see that $\ell=3$ or $4$. If $\ell=3$, then $v=(t_1, t_2, 1)$ where $\lg t_1,t_2\rg\cong\mz_2\times\mz_2$, and then $S=\{(t_1,t_2,1)^{\bf t}\mid {\bf t}=(t,t,t)\in N_u\}$ and $\G\cong\G_{U(3,4)^3}$ (see Example~\ref{exam:sd-u34-1}). If $\ell=4$, then $v=(t_1, t_2, t_3, 1)$ where $\lg t_1,t_2,t_3\rg\cong\mz_2\times\mz_2$, and then $S=\{(t_1,t_2,t_3, 1)^{\bf t}\mid {\bf t}=(t,t,t,t)\in N_u\}$ and $\G\cong\G_{U(3,4)^4}$ (see Example~\ref{exam:sd-u34-2}).\hfill\qed

\section{Examples of quasiprimitive $3$-CH graphs}\label{sec:examples}
In this section, we shall show that there are examples of finite $(G,3)$-CH graphs such that $G$ is quasiprimitive on the vertices of type HA, AS, PA or TW, and hence prove part (3) of Theorem~\ref{th-quasiprimitive}. Throughout, $\G$ will be neither complete nor complete multipartite and have girth $3$.

We first consider the case where $G$ is of HA-type. From Lemma~\ref{normal cayley} we can obtain the following result.

\begin{lem}\label{lem-HA-disconnected}
Suppose that $\G=(V,E)$ is a $(G,3)$-CH graph such that $G$ is quasiprimitive on $V$ of type {\rm HA}. Then $\G$ is a Cayley graph, say $\Cay(H,S)$, on the socle $H$ of $G$. Furthermore, one of the following holds.
\begin{enumerate}
  \item [{\rm (1)}]\ $[S]\cong r\K_2$, $H\cong\mz_3^{r_1} (2\leq r_1\leq r)$.
  \item [{\rm (2)}]\ $[S]\cong r\K_{2^n-1}$, $H\cong\mz_2^{r_1n} (2\leq r_1\leq r)$.
  \item [{\rm (3)}]\ $[S]$ is connected, $H\cong\mz_2^n$ for some integer $n$.
\end{enumerate}
\end{lem}




The following three examples show that all the three possible cases in Lemma~\ref{lem-HA-disconnected} do occur.

\begin{exam}\label{HA-type-disconnected3}$(${\rm HA}, locally disconnected$)$\

{\rm Let $H_i\cong\mz_2^\ell$ with $i=1,\ldots,r$ and $\ell\geq 2$, and let $H=H_1\times\cdots\times H_1$ and $S=\bigcup_{i}^r (H_i-\{1_{H_i}\})$. Define $\G=\Cay(H,S)$. Then $\Aut(H,S)\cong \GL(\ell,2)\wr S_r$, and $\G$ is $(G, 3)$-CH, where $G=H_R\:\Aut(H, S)$.
Suppose that there is a subgroup, say $N$, of $H$ which is normalised by $\Aut(H,S)$. Let $H_i=\lg a_{i1}\rg\times\cdots\times\lg a_{i\ell}\rg$ for $i=1,\ldots, r$. Take a non-identity element $g$ in $N$ so that $g=\prod_{i=1}^r a_{i1}^{x_{i1}}a_{i2}^{x_{i2}}\cdots a_{i\ell}^{x_{i\ell}}$ with $x_{ij}\in\mz_2$. Without loss of generality, assume that $a_{11}^{x_{11}}a_{12}^{x_{12}}\cdots a_{1\ell}^{x_{1\ell}}\not=1$. We may choose an $\a\in\Aut(H, S)=\GL(\ell,2)\wr S_r$ such that $H_1^\a=H_1$ and $(a_{11}^{x_{11}}a_{12}^{x_{12}}\cdots a_{1\ell}^{x_{1\ell}})^\a\neq a_{11}^{x_{11}}a_{12}^{x_{12}}\cdots a_{1\ell}^{x_{1\ell}}$ while $(a_{ij}^{x_{ij}})^\a=a_{ij}^{x_{ij}}$ for all $i>1$. So, we have $1\neq gg^\a\in H_1-\{1_{H_1}\}\subseteq S$. Since $\Aut(H,S)$ is transitive on $S$, one has $S\subseteq N$. It follows that $N=\lg S\rg=H$, and so $G$ is primitive on the vertex set $V\G$ of $\G$. It is also easy to see that the subgraph induced by $S$ is isomorphic to $r\K_{2^\ell-1}$.}
\end{exam}

\begin{exam}\label{HA-type-disconnected2}$(${\rm HA}, locally disconnected$)$\

{\rm Let $H=\lg a_1\rg\times\cdots\times\lg a_r\rg\cong\mz_3^{r}$, and let $S=\{a_i,a_i^{-1}\ |\ i=1,\ldots,r\}$. Define $\G=\Cay(H,S)$. Then $\Aut(H,S)\cong S_{2}\wr S_r$, and $\G$ is $(G, 3)$-CH, where $G=H_R\:\Aut(H, S)$. Suppose that there is a subgroup, say $N$, of $H$ which is normalised by $\Aut(H,S)$. Take a non-identity element $g$ in $N$ so that $g=a_1^{i_1}a_2^{i_2}\cdots a_{r}^{i_r}$ for some $i_j\in\mz_3$. Without loss of generality, assume that $i_1\neq 0$. We can take an $\a\in\Aut(H,S)\cong S_{2}\wr S_r$ such that $a_1^\a=a_1$ and $a_i^\a=a_i^{-1}$ with $2\leq i\leq r$. So, $g^\a=a_1^{i_1}a_2^{-i_2}\cdots a_{r}^{-i_r}\in N$, and hence $1\neq gg^\a=a_1^{2i_1}\in N$. Since $\Aut(H, S)$ is transitive on $S$, all $a_i$'s are in $N$, and so $N=H$. Thus, $H$ is a minimal normal subgroup of $G$, and so $G$ is primitive on the vertex set $V\G$ of $\G$. It is also easy to see that the subgraph induced by $S$ is isomorphic to $r\K_{2}$.}
\end{exam}

\begin{exam}\label{HA-type-connected}$(${\rm HA}, locally connected$)$\

{\rm Let $H=\lg a_1\rg\times\cdots\times\lg a_{n-1}\rg\cong\mz_2^{n-1}$ and $\Om=\{a_1, \ldots, a_{n-1}, a_1a_2\cdots a_{n-1}\}$, where $n\geq 3$ is odd. Then the Cayley graph $\Cay(H,\Om)$ is the {\em folded cube}, the antipodal quotient of the $k$-cube graph. Set $G=H_R\:\Aut(H, \Omega)$. By \cite{Ivanov-Praeger}, $\Aut(H,\Om)\cong S_n$ and $G$ is primitive on $V\Cay(H,\Om)$. Let $S$ be the set of products of $(n-2)$ elements in $\Om$. Clearly, $S$ is an orbit of $\Aut(H,\Om)$. In fact, $\Aut(H,\Om)$ acts primitively on $S$ with rank $3$ (see \cite{Bannai}). Define $\G=\Cay(H, S)$. Then $\G$ is an orbital graph of $G$ on $H\times H$. So, $\G$ is connected. It is easy to see that $\G$ has girth $3$.
By Proposition~\ref{prop-local-action}, $\G$ is $(G, 3)$-CH, and by Theorem~\ref{local-pty}, $\G$ is locally connected. }
\end{exam}

We now look at the case where $G$ is of AS-type. First, we provide a family of locally connected graphs.

\begin{exam}\label{AS-type-connected}$(${\rm AS}, locally connected$)$\

{\rm The {\em Johnson graphs} $J(n,k)$ have, as vertex set, the set $V$ of $k$-element subsets of the set $\{1,2,\ldots,n\}$, for some $k, 1\leq k\leq n-1$, and two $k$-elements subsets $\a,\b$ are adjacent if and only if the intersection $\a\cap\b$ has size $k-1$. The valency of $J(n,k)$ is $k(n-k)$ and $\Aut(J(n,k))=S_n$. (See \cite[Section~9.1]{Brouwer-Book}.)

Let $n\geq7$. Consider the complementary graph $J(n,2)^c$ of $J(n,2)$ and let $G=\Aut(J(n,2))=S_n$. Then $G$ is primitive on $V$ with rank $3$ (see \cite{Bannai}). For $\a=\{1,2\}\in V$, we have $G_\a=S_{n-2}\times S_2$. It is easy to see the neighbourhood $N$ of $\a$ in $J(n,2)^c$ is the set of $2$-element subsets of the set $\{3,4,\ldots,n\}$. Again, by \cite{Bannai}, $G_\a$ is primitive on $N$ with rank $3$. It is easy to see that $J(n,2)^c$ has girth $3$. By Proposition~\ref{prop-local-action}, $\G$ is $(G, 3)$-CH, and by Theorem~\ref{local-pty}, $\G$ is locally connected. }
\end{exam}

To construct the locally disconnected examples with $G$ is of AS-type, we need the following lemma.

\begin{lem}\label{lem-line-graph}
Let $\G'$ be a non-complete $(G,3)$-arc-transitive graph such that $G_u$ is $3$-transitive on $\G'(u)$ for any $u\in V\G'$. Let $\G$ be the line graph of $\G'$. Then $\G$ is a $(G,3)$-CH graph.
\end{lem}

\f\demo Clearly, $\G=E\G'$. Take $e=\{u,v\}\in E\G'$. Then $[\G_1(e)]\cong 2\K_{k-1}$, where $k$ is the valency of $\G'$. Since $G_u$ is $3$-transitive on $\G'(u)$, $G_{uv}$ is $2$-transitive on both $\G'(u)-\{v\}$ and $\G'(v)-\{u\}$. For any $w\in\G'(u)-\{v\}$, $G_{uvw}$ is
transitive on $\G'(u)-\{v\}$. Since $\G$ is $(G,3)$-arc-transitive, $G_{uvw}$ is also transitive on $\G'(v)-\{u\}$.
Thus, $G_e^{\G_1(e)}$ is a transitive permutation group of rank $3$. By Proposition~\ref{prop-local-action}, $\G$ is a $(G,3)$-CH graph.\hfill\qed

\begin{exam}\label{AS-type-disconnected}$(${\rm AS}, locally disconnected$)$\

{\rm Let $\G'=J(5,2)^c$ be the Petersen graph and $G=\Aut(\G')=S_5$. Then for any $u\in V\G'$, $G_u$  is $3$-transitive on $\G'(u)$. By Lemma~\ref{lem-line-graph}, the line graph of $\G'$ is a $(G,3)$-CH graph.}
\end{exam}

Examples where $G$ is of type PA also occur.

\begin{exam}\label{PA-type-connected}$(${\rm PA}, locally connected, locally disconnected$)$

{\rm Let $\D=\{0,1,\ldots,k-1\}$ and let $d\geq 2$ be an integer. The {\em Hamming graph} $H(d,k)$ has vertex set $\D^d$, the set of ordered $d$-tuples of elements of $\D$, or sequences of length $d$ from $\D$. Two vertices are adjacent if they differ in precisely one coordinate. The Hamming graph has valency $d(k-1)$ and diameter $d$. It is $G$-distance transitive with $G=\Aut(H(d,k))=S_k\wr S_d$, and it is $G$-vertex-primitive (of type PA) if and only if $k\geq 3$ (see \cite[Section~9.2]{Brouwer-Book}). Let $k\geq 4$, $d=2$ and $\G=H(2,k)$. Set $\a=(0,0)$. Then $G_\a=S_{k-1}\wr S_2$ which is transitive on $\G_1(\a)$ with rank $3$ and primitive on $\G_2(\a)$ with rank $3$. Clearly, $[\G_1(\a)]\cong 2\K_{k-1}$. By Proposition~\ref{prop-local-action}, both $\G$ and its complementary graph $\G^c$ are $(G,3)$-CH. Furthermore, $\G$ is locally disconnected while $\G^c$ is locally connected.}
\end{exam}

Finally, we give an example in case $G$ is of type TW.

\begin{exam}\label{TW-type-disconnected}
$(${\rm TW} locally disconnected$)$ 

{\rm Let $h_1, h_2, a, p, q\in { S_{14}}$ be such that
\[\begin{array}{l}
h_1= (1,2,3), h_2 = (1,4,5)(2,6,7),\\

a=(1,8)(2,9)(3,10)(4,11)(5,12)(6,13)(7,14),\\

p= (4, 5)(6, 7)(8, 9)(11, 13)(12, 14),\\

q= (1, 2)(4, 6)(5, 7)(11, 12)(13, 14).
\end{array}
\]
Let $h= h_1h_2^a$, $G=\lg h, a, p, q\rg$ and $H=\lg h, p, q\rg$. By {Magma}~\cite{BCP}, we have $G\cong { S_7}\wr { S_2}$, $H\cong S_3\times S_2$, $|H: H^a\cap H|=3$. Let $\G'$ be the graph with vertex set the set of right cosets of $H$ in $G$, and two different cosets $Hg_1, Hg_2$ are adjacent if and only if $g_2g_1^{-1}\in HaH$. Then $\G'$ has valency $3$ and the stabiliser $G_u\cong S_3\times S_2$, where $u=H$. So $\G'$ is $3$-arc-transitive and $G_u$ is $3$-transitive on $\G'(u)$. By Lemma~\ref{lem-line-graph}, the line graph $\G$ of $\G'$ is $(G,3)$-CH. Clearly, $\soc(G)=A_7\times A_7$ and the edge stabiliser $G_{\{u,v\}}=\lg p,q\rg:\lg a\rg\cong D_8$, where $v=Ha$. Furthermore, $\soc(G)\cap G_{\{u,v\}}=1$. It follows that $G$ is regular on the edges of $\G'$ and so regular on the vertices of $\G$. This implies that $G$ is quaisprimitive on $V\G$ of type TW.}
\end{exam}

\section{Proof of Theorem~\ref{th-quasiprimitive}}\label{sec:proof}

Let $\G=(V,E)$ be a $(G,3)$-CH graph which is neither complete nor complete multipartite, and $G$ is quasiprimitive on $V$. By Theorem~\ref{th:hs-hc-cd}, $G$ is not of type {\rm HS, HC} or {\rm CD}, and so part (1) of Theorem~\ref{th-quasiprimitive} holds. By Theorem~\ref{SD}, if $G$ is quasiprimitive on $V$ of type SD, then $\G\cong\G_{A_5^3}, \G_{A_5^4}, \G_{U(3,4)^3}$ or $\G_{U(3,4)^4}$. This proves part (2) of Theorem~\ref{th-quasiprimitive}. From Examples~\ref{HA-type-disconnected3}--\ref{TW-type-disconnected}, there do exist Examples of $(G,3)$-CH graphs in case $G$ is of type {\rm HA, AS, PA} or {\rm TW}. This implies part (3) of Theorem~\ref{th-quasiprimitive}.

\medskip
\f {\bf Acknowledgements:} This work was supported by the National
Natural Science Foundation of China (12071023,1211101360).

\end{document}